\documentclass[fleqn,11pt,twoside]{article}

\usepackage{amsthm,amsthm,amssymb, color, xcolor,epsfig, graphics, subfigure,cite}

\usepackage{amsmath, graphicx, latexsym, lscape }

\makeatletter
\newcommand{\copyrightnote}[2]{{\renewcommand{\thefootnote}{}
 \footnotetext{\small\it
\begin{flushleft}
 \copyright \ #1   #2  
\end{flushleft}}}}

\newcommand{\Name}[1]{\begin{flushleft}
                       \LARGE \bf #1
                       \end{flushleft}\vspace{-3mm}}

\newcommand{\Author}[1]{\begin{flushleft}
                       \it #1 \end{flushleft}}

\newcommand{\Address}[1]{\begin{flushleft}
                       \it #1 \end{flushleft}}

\newcommand{\Date}[1]{\begin{flushleft}
                      \small  \it #1 \end{flushleft}}

\newcommand{\evenhead}{Author \ name}
\newcommand{\oddhead}{Article \ name}

\renewcommand{\@evenhead}{
\hspace*{-3pt}\raisebox{-15pt}[\headheight][0pt]{\vbox{\hbox to \textwidth
{\thepage \hfil \evenhead}\vskip4pt \hrule}}}
\renewcommand{\@oddhead}{
\hspace*{-3pt}\raisebox{-15pt}[\headheight][0pt]{\vbox{\hbox to \textwidth
{\oddhead \hfil \thepage}\vskip4pt\hrule}}}
\renewcommand{\@evenfoot}{}
\renewcommand{\@oddfoot}{}

\long\def\@makecaption#1#2{%
  \vskip\abovecaptionskip
  \sbox\@tempboxa{\small \textbf{#1.}\ \ #2}%
  \ifdim \wd\@tempboxa >\hsize
    {\small \textbf{#1.}\ \ #2}\par
  \else
    \global \@minipagefalse
    \hb@xt@\hsize{\hfil\box\@tempboxa\hfil}%
  \fi
  \vskip\belowcaptionskip}

\newcommand{\JNMPnumberwithin}[3][\arabic]{%
  \@ifundefined{c@#2}{\@nocounterr{#2}}{%
    \@ifundefined{c@#3}{\@nocnterr{#3}}{%
      \@addtoreset{#2}{#3}%
      \@xp\xdef\csname the#2\endcsname{%
        \@xp\@nx\csname the#3\endcsname .\@nx#1{#2}}}}%
}

\newcommand{\resetfootnoterule} {
  \renewcommand\footnoterule{%
  \kern-3\p@
  \hrule\@width.4\columnwidth
  \kern2.6\p@}
}

\renewcommand{\footnoterule}{}

\makeatother

\theoremstyle{definition}

\newcommand{\blu}{\color{blue}}

\def\gcite#1{\cite{#1}}

\def\a{\alpha}
\def\b{\beta}
\def\ga{\gamma}
\def\de{\delta}   
\def\eps{\varepsilon}
\def\vphi{\varphi}
\def\la{\lambda}

\def\s{\sigma}

\def\om{\omega}

\def\vphi{\varphi}

\def\F{{\cal F}}
\def\G{{\cal G}}

\def\S{{\mathcal S}}
\def\T{{\rm T}}
\def\V{{\cal V}}

\def\De{\Delta}

\def\pa{\partial}

\def\ss{\subset}
\def\sse{\subseteq}

\def\<{\langle}
\def\>{\rangle}

\def\xd{\dot{x}}
\def\yd{\dot{y}}
\def\zd{\dot{z}}
\def\xb{{\bf x}}
\def\yb{{\bf y}}
\def\ub{{\bf u}}

\def\({\left(}
\def\){\right)}
\def\[{\left[}
\def\]{\right]}
\def\=#1{\bar #1}
\def\~#1{\widetilde #1}
\def\wt#1{\widetilde #1}
\def\.#1{\dot #1}
\def\^#1{\widehat #1}
\def\wh#1{\widehat #1}
\def\"#1{\ddot #1}

\def\bar#1{\overline{#1}}

\def\eeq{\end{equation}}
\def\beq{\begin{equation}}

\def\beql#1{\begin{equation} \label{#1}}

\def\EOR{\hfill $\odot$}
\def\EOP{\hfill $\diamondsuit$}
\def\EOE{\hfill $\Box$}

\def\eqref#1{(\ref{#1})}

\def\symref{AVL,CGb,KrV,Olv1,Olv2,Ovs,Ste,Win1,Win2}

\begin{document}

\renewcommand{\evenhead}{ {\LARGE\textcolor{blue!10!black!40!green}{{\sf \ \ \ ]ocnmp[}}}\strut\hfill G Gaeta}
\renewcommand{\oddhead}{ {\LARGE\textcolor{blue!10!black!40!green}{{\sf ]ocnmp[}}}\ \ \ \ \   Conditional symmetries for dynamical systems}

\thispagestyle{empty}
\newcommand{\FistPageHead}[3]{
\begin{flushleft}
\raisebox{8mm}[0pt][0pt]
{\footnotesize \sf
\parbox{150mm}{{Open Communications in Nonlinear Mathematical Physics}\ \  \ {\LARGE\textcolor{blue!10!black!40!green}{]ocnmp[}}
\ \ Vol.5 (2025) pp
#2\hfill {\sc #3}}}\vspace{-13mm}
\end{flushleft}}

\FistPageHead{1}{\pageref{firstpage}--\pageref{lastpage}}{ \ \ Article}

\strut\hfill

\strut\hfill

\copyrightnote{The author(s). Distributed under a Creative Commons Attribution 4.0 International License}

\Name{Conditional symmetries and conditional constants of motion for dynamical systems}

\Author{Giuseppe Gaeta$^{1,2}$}

\Address{${}^1$  Dipartimento di Matematica, Universit\`a degli Studi di Milano, via Saldini 50, 20133 Milano (Italy) {\rm and} INFN Sezione di Milano; {\tt giuseppe.gaeta@unimi.it} \\ 
${}^2$ {\it SMRI, 00058 Santa Marinella (Italy)}}

\Date{Received December 4, 2024; Accepted January 6, 2025}

\setcounter{equation}{0}

\begin{abstract}
\noindent 
Conditional symmetries were introduced by Levi and Winternitz in their 1989 seminal paper to deal with nonlinear PDEs. Here we discuss their application in the framework of ODEs, and more specifically Dynamical Systems; it turns out they are closely related to two established -- albeit maybe less widely known -- concepts, i.e. orbital symmetries and configurational invariants. The paper is devoted to studying the interplay of these notions, and their application in the study of Dynamical Systems, with special attention to invariant manifolds of these. 

\bigskip\noindent\emph{{\blu This paper is dedicated to the dear memory of Decio Levi and Pavel Winternitz}}
\end{abstract}

\label{firstpage}




\addcontentsline{toc}{section}{\ \ \ \ Introduction}

\section*{Introduction}

The purpose of this work is to reconsider some notions, all of them related to symmetries of differential equations and dynamical systems, but which appear so far to have been considered each by itself, and to investigate their mutual relations.

These notions are the well known one of \emph{conditional symmetries} introduced by Levi and Winternitz \gcite{LeWin,LWR} and mainly applied to PDEs (and to partial difference equations, a theme we will not consider here); the notion of \emph{orbital symmetries} introduced by Walcher \gcite{WalOS,WalMPS} in the study of dynamical systems; and the notion of \emph{configurational invariants} -- also called {\it conditional constants of motion} -- introduced by Sarlet, Leach and Cantrijn \gcite{SLC} in their study of ``weak integrability'' (or rather, in the nomenclature we will adopt here, conditional integrability) of dynamical systems.

More specifically, our main theme will be to consider how the notion of conditional symmetries applies to \emph{dynamical systems} (which seems not to have undergone specific investigation so far); we will then show that the notions of orbital symmetries (defined for dynamical systems) and configurational invariants (defined for ODEs) recalled above  have a rich interplay among themselves and with conditional symmetries.

In Section \ref{sec:disc} we will discuss several possible extensions (and, at the same time, limitations) of our discussion.

Albeit we cannot claim any substantial new result, we trust that a better understanding of the interplay between these notions will help to profitably apply them in the study of Dynamical Systems.

\bigskip

Summation over repeated indices will be routinely understood. By ``smooth'' we will always understand $\mathcal{C}^\infty$. The end of Examples, Proofs and Remarks will be signalled by the symbols $\Box$, $\diamondsuit$ and $\odot$ respectively.

\section{Lie-point symmetries for dynamical systems}

We start by recalling some basic facts about symmetry analysis of differential equations \gcite{\symref}; as we are mostly interested in dynamical system \gcite{ArnODE,ArnEMS,GuckHolm,Ver1}, we will discuss the general theory in this context -- which allows for a slightly lighter notation than in the general case.

Let $M$ be a smooth manifold, which we think of as embedded in $R^N$, the latter being equipped with coordinates $(x^1,...,x^N)$. In the following we also write
$$ \pa_i \ := \ \frac{\pa}{\pa x^i} \ ; $$
the vector operator of components $(\pa_1,...,\pa_n)$ will be denoted as $\nabla$.

A \emph{dynamical system} on $M$ is a system of first order autonomous ODEs
\beql{eq:DS} \xd \ = \ f (x) \eeq
where $f: M \to \T M $ is a smooth (i.e. $\mathcal{C}^\infty$) function.

\subsection{Symmetry algebra}
\label{sec:s1}

The dynamical system \eqref{eq:DS} corresponds to the flow on $M$ under the action of the vector field 
\beql{eq:13} \wh{f} \ = \ f(x) \, \nabla \ = \ f^i (x) \ \frac{\pa}{\pa x^i} \ . \eeq

The set of vector fields on $M$ equipped with the Lie commutator $[.,.]$ form a Lie algebra, which we denote by $\mathcal{V}$. If $\^f$ is as in \eqref{eq:13} and we consider a new vector field
\beql{eq:sVF} \^s \ = \ s^i (x) \, \pa_i \ , \eeq their commutation is given by
\beql{eq:comm} [\^s,\^f] \ = \ \( \^s (f^i) \ - \ \^f (s^i ) \) \, \pa_i \ = \ \( s^j \pa_j f^i \ - \ f^j \pa_j s^i \) \, \pa_i \ . \eeq

Clearly, to any vector field $\wh{f} \in \V$ we can associate a smooth vector function $f : M \to \T M$, i.e. a section of the \emph{tangent bundle} $\T M$ for $M$. The space of such smooth functions will be denoted by $\F$; this is naturally equipped with a Lie-Poisson bracket\footnote{Needless to say, this is not to be mistaken with the Poisson bracket met in Hamiltonian Dynamics.} $\{ . , . \} : \F \times \F \to \F$ induced by $[.,.]$; this is given by
\beql{eq:14} \{ f , g \} \ := \ (f \cdot \nabla) \, g \ - \ (g \cdot \nabla ) \, f \ . \eeq

We can now define symmetry vector fields. It should be mentioned that what we call a ``symmetry vector field'', see Definition 1, is, more precisely, the generator of a one-parameter local group of symmetries; we will adopt this  abuse of language (rather common in the literature) for ease of discussion.

\medskip\noindent
{\bf Definition 1.}
{\it A vector field $\^s \in \V$ is a \emph{Lie-point time-independent} (LPTI) symmetry of the dynamical system \eqref{eq:DS} if under the flow of $\^s$ solutions to \eqref{eq:DS} are transformed into (generally, different) solutions to \eqref{eq:DS}. }

\medskip\noindent
{\bf Lemma 1.} {\it In the present setting, $\^s$ is a symmetry for \eqref{eq:DS} if and only if
\beql{eq:15c} [ \^f , \^s ] \  = \ 0 \ ; \eeq
or, in terms of \eqref{eq:14},
\beql{eq:15} \{ f , s \} \ = \ 0 \ . \eeq}

\medskip\noindent
{\bf Proof.} See any book on symmetries of differential equations \gcite{\symref}. \EOP
\bigskip

It is easy to see that the vector fields satisfying \eqref{eq:15c} form a \emph{Lie algebra}; the same holds for functions satisfying \eqref{eq:15} or even \eqref{eq:15p}. We will denote these, respectively, by $\^\G_{\^f} \ss \V$ and by $\G_f \ss \F$.

\medskip\noindent
{\bf Remark 1.} Here we are considering not only time-autonomous dynamical systems, but also (symmetry) vector fields which do not act nor depend on time, i.e. Lie-Point Time-Independent (LPTI) symmetries. The case of general (i.e. non-autonomous and acting on time) vector fields -- in which we deal with $\overline{M} = R \times M$ and with $\overline{\V}$ (the set of vector fields on $\overline{M}$) and $\overline{\F}$ (the set of smooth functions from $\overline{M}$ to $\T \overline{M}$) -- is readily recovered by adding a new coordinate $x^0$ whose evolution is given by $\xd^0 = 1$. If $f$ is autonomous but $s$ is general, i.e.
\beql{eq:sVFgen} \^s \ = \ s^i (x,t)  \, \pa_i \ , \eeq
the relation in terms of commutators of vector fields is just the same as \eqref{eq:15c}, while \eqref{eq:15} is replaced by
\beql{eq:15p} s_t \ + \ \{f , s \} \ = \ 0 \ . \eeq

More generally \gcite{\symref}, LP symmetries would be vector fields on $R \times M$ of the form
\beql{eq:16} \eta \ = \ \tau (x,t) \, \frac{\pa}{\pa t} \ + \ \vphi^i (x,t) \, \frac{\pa}{\pa x^i} \ . \eeq The determining equations for symmetries are then (with summation over repeated indices)
\beq \pa_t \vphi^i \ + \ \( f^j \pa_j \vphi^i \ - \ \vphi^j \pa_j f^i \) \ = \ f^i \ \( \pa_t \tau \ + \ f^j\pa_j \tau \) \ . \eeq
For $\tau = c$, corresponding to time translations (and including the case $\tau = 0$ of interest here), we just get
\beq \pa_t \, \vphi^i \ + \ \{ f , \vphi \}^i \ = \ 0 \ , \eeq
which is just \eqref{eq:15p}. \EOR

\medskip\noindent
{\bf Remark 2.} If $\tau$ is nowhere zero, we can divide \eqref{eq:16} by $\tau$, reaching the form $\pa_t + \vphi \nabla$; if $\tau$ has some zero, proceeding in this way produces a vector field which is singular at such points. \EOR

\medskip\noindent
{\bf Remark 3.} In the study of (autonomous, and not only) dynamical systems it is specially convenient to consider vector fields not acting on time, i.e. with $\tau \equiv 0$. One reason for this is the following: we know that a vector field is a symmetry for a given equation if and only if its evolutionary representative is \gcite{Olv1}; for a vector field of the form \eqref{eq:sVF}, the evolutionary representative is
\beql{eq:sVFev} \^s_e \ = \ Q^i \, \pa_i \ , \ \ \ \ Q^i \ = \ s^i \ - \ \tau \, \xd^i \ . \eeq When we restrict this to the solution manifold for \eqref{eq:DS} we just get
\beq Q^i \ = \ s^i \ - \ \tau \, f^i \ . \eeq Thus, as $s^i (x,t)$ are arbitrary functions, the term $\tau (x,t) f^i (x)$ can be absorbed into the choice of $s^i (x,t)$. Note that the same argument holds for non-autonomous dynamical systems as well. For a more general discussion see e.g. \gcite{Cds,CGihp,CGncb,G93}. \EOR
\bigskip

\subsection{Module structure of the symmetry algebra}
\label{sec:s2}

We now consider the case where \eqref{eq:DS} -- or $f$ for short --  admits a \emph{constant of motion} (CM); that is, there is a smooth scalar function $P(x): M \to R$ which is constant under the flow of \eqref{eq:DS}, i.e. such that
\beql{eq:21} \frac{d P}{dt} \ = \ (\xd \cdot \nabla) \, P \ = \ \( f \cdot \nabla \) \, P \ = \ 0 \ . \eeq

Note that for a given dynamical system these are determined by solving the associate characteristic system
$$ \frac{d x^1}{f^1} \ = \ ... \ = \ \frac{d x^n}{f^n} \ . $$

It turns out that in this case the set of symmetries for the dynamical system has, beside the structure of Lie algebra, also that of Lie module .(For definition and properties of modules, see any text in Algebra; or e.g. \gcite{Kir}.) In other words, we have:

\medskip\noindent
{\bf Lemma 2.} {\it Let the dynamical system \eqref{eq:DS} admit $P$ as a CM. If $\^s$ is a symmetry for $\^f$, then $\^r = g(P) \^s$ is also a symmetry for $\^f$, for any smooth scalar function $g: R \to R$.}

\medskip\noindent
{\bf Proof.} In fact,
\beql{eq:22} \{f,g(P)s\} \ = \ g'(P) \, s \, ( f \cdot \nabla) P \ + \ g(P) \, \{ f , s \} \ = \ g(P) \, \{ f , s \} \ , \eeq
where we have used first the assumption that $P$ is a CM, and then $\{f,s\}=0$. This concludes the proof for LPTI symmetries. As for time-dependent ones, we note that the time derivative is just
\beq \pa_t r^i \ = \ \pa_t \( g(P) \, s^i \) \ = \ g(P) \, \pa_t s^i \ + \ \( \pa_t g(P)\) \, s^i \ ; \eeq
on the flow on $\^f$, $\pa_t g(P) = 0$ by the assumption $P$ is a CM, and thus eq. \eqref{eq:15p} reads
\beq g(P) \ \( s^i_t \ + \ f^j \pa_j s^i \ - \ s^j \pa_j f^i \) \ ; \eeq this vanishes if the term in brackets does.  \EOP
\bigskip

The CM for a given dynamical system form an associative and commutative algebra, closed under sum and product: if $P_1$ and $P_2$ are CM for $f$, obviously also $\wt{P} = P_1 + P_2$ is a CM for $f$, and similarly for $\wh{P} = P_1 P_2$. We denote this algebra by $I_f$. Note that the set $I^{(0)}$ of functions which are constant on $M$ is trivially in $I_f$ for any $f$; we call $I_f/I^{(0)}$ the nontrivial algebra of CM for $f$.

\medskip\noindent
{\bf Lemma 3.} {\it The set $\G_f$ is a Lie module over the algebra $I_f$.}

\medskip\noindent
{\bf Proof.} We have to show that if $s_1,s_2 \in \G_f$ and $P_1,P_2 \in I_f$, then $\s = P_1 s_1 + P_2 s_2 \in \G_f$ as well. This follows from the bilinearity of the Lie-Poisson bracket $\{.,.\}$ via a trivial computation.

\medskip\noindent
{\bf Corollary.} {\it If $I_f$ is nontrivial, then $\G_f$ is infinite dimensional as a Lie algebra.}

\medskip\noindent
{\bf Remark 4.} Note that in this case $\G_f$ will be in general (for $f \not\equiv 0$) finite dimensional as a module. \EOR

\subsection{Invariant solutions}
\label{sec:IF}

In the following we will deal with invariant solutions, i.e. functions which are solutions of a differential equations and which are invariant under the action of a vector field. We recall that in full generality, i.e. when we consider functions $u^a (x^1,...,x^n)$ depending on several independent variables and correspondingly vector fields
\beq \^v \ = \ \vphi^a (x,u) \, \frac{\pa}{\pa u^a} \ + \ \xi^i (x,u) \, \frac{\pa}{\pa x^i} \ , \eeq
the function $u^a = f^a (x)$ is invariant under $\^v$ if and only if the equations
\beq \vphi^a \ = \ \xi^i \, f^a_i \eeq are satisfied.

This descends immediately from the fact that under the infinitesimal action of $\^v$ (i.e. under $\exp[ \eps \^v ]$) the function $f$ is mapped into $\wt{f}$ with
\beq \wt{f}^a \ = \ f^a \ + \ \eps \ \[ \vphi^a \ - \ \xi^i \, u^a_i \] \ . \eeq
We refer e.g. to \gcite{\symref} for details.

\subsection{Orbital symmetries}
\label{sec:OS}

In the case of dynamical systems, it makes full sense to consider, beside ordinary symmetries, also so called \emph{orbital symmetries} \gcite{WalOS,WalMPS}. While a symmetry preserves the equation, an orbital symmetry preserves the equation up to a time reparametrization. In the case of dynamical systems, this means that the equation can be transformed into a different one whose solutions have the same \emph{trajectories} as the original equation. We refer to \gcite{WalOS,WalMPS} for a general discussion of orbital symmetries and their applications. (As orbital symmetries are likely to be less widely known than standard symmetries, some simple examples of the former are collected in Appendix A to help the reader not familiar with these to fix ideas.)

\medskip\noindent
{\bf Definition 2.} {\it Given the dynamical system \eqref{eq:DS}, the vector field \eqref{eq:sVF} or more generally a vector field \eqref{eq:16} is an \emph{orbital symmetry}  for \eqref{eq:DS} if it maps solutions' trajectories into (generally, different) solutions' trajectories.}
\bigskip

Note that vector fields inducing just a reparametrization of time are trivial orbital symmetries; similarly and for the same reason, in any vector field of the general form \eqref{eq:16} we can drop the part along $\pa_t$, i.e. we can set $\tau = 0$, obtaining a vector field which is fully equivalent for what concerns being an orbital symmetry (for any dynamical system). We will thus concentrate only on nontrivial orbital symmetries.

\medskip\noindent
{\bf Lemma 4.} {\it The orbital symmetries for \eqref{eq:DS} are the vector fields \eqref{eq:sVFgen} which satisfy
\beql{eq:23} \[ \^s , \^f \] \ = \ \la \, \^f \ ; \eeq
or, in terms of components,
\beql{eq:24} s^i_t \ + \ \{ s , f \}^i \ = \ \la \, f^i \ . \eeq }

\medskip\noindent
{\bf Proof.} Given the dynamical system \eqref{eq:DS}, the dynamical systems having the same trajectories are those which are written as
\beql{eq:22c} \xd^i \ = \ g^i(x,t) \ = \ \la (x,t) \ f^i (x,t) \ , \ \ \ \la (x,t) \not= 0 \ \forall x,t \ ; \eeq note that $\la (x,t)$ is a (nowhere zero) scalar function; this is sometimes also written as
$$ \la (x,t) \ = \ \exp [ \mu (x,t) ] \ , $$ in order to emphasize it is nowhere vanishing.
A direct computation shows that this is equivalent to the condition given in the statement. \EOP

\medskip\noindent
{\bf Lemma 5.} {\it The set of orbital symmetries for a given dynamical system has the structure of a Lie algebra, and also that of a Lie module over the ring of CM for the dynamical system.}

\medskip\noindent
{\bf Proof.} We prove these statements in the vector fields language; the proof in terms of component functions is analogous. We also only consider nontrivial orbital symmetries, the extension being immediate.

If $\^s$ and $\^r$ are both orbital symmetries, so that there are nonzero functions $\s$ and $\rho$ such that
$$ [ \^s,\^f] \ = \ \s \, \^f \ , \ \ \ [\^r , \^f ] \ = \ \rho \, \^f \ , $$
the using Jacobi identity we obtain
$$ \[ [\^s,\^r] , \^f \] \ = \ \( \^s (\rho) \ - \ \^r (\s ) \) \, \^ f \ , $$ thus showing that $[\^s,\^r]$ is indeed again an orbital symmetry.

Now, let $\mu : R^n \to R$ be a CM and $\^s$ an orbital symmetry for the dynamical system \eqref{eq:DS}; then
$$ [ \mu \^s , \^f] \ = \ \mu \, [\^s,\^f] \ - \ \( \^f (\mu ) \) \, \^s \ = \ \mu \, \la \, \^f \ , $$ thus showing that $\mu \^s$ is again an orbital symmetry.  \EOP

\medskip\noindent
{\bf Lemma 6.} {\it For the dynamical system $\xd^i = f^i (x)$, any vector field of the form $\^v = \vphi^i (x,t) \pa_i$ with
\beq \vphi^i (x,t) \ = \ \theta (x,t) \ f^i (x) \eeq with $\theta$ a smooth nowhere zero function is an orbital symmetry.}

\medskip\noindent
{\bf Proof.} By a direct computation,
\begin{eqnarray*}
\vphi^i_t \ + \ \{f, \vphi\}^i &=& \vphi^i_t \ + \ f^j \, \pa_j \vphi^i \ - \ \vphi^j \, \pa_j f^i \\
&=& \theta_t \, f^i \ + \ f^j \, (\theta_j f^i \, + \, \theta f^i_j ) \ - \ \theta \, f^j \, f^i_j \\
&=& \( \theta_t \ + \ f^j \, \theta_j \) \ f^i \ . \end{eqnarray*}
Note that the function $\la (x,t)$ appearing in the definition \eqref{eq:23}, \eqref{eq:24} of orbital symmetries is in this case just the total time derivative of $\theta (x,t)$ computed on the flow of the dynamical system, i.e.
$$ \la (x,t) \ = \ D_t \, \theta \ = \ \theta_t \ + \ f^j \, (\pa_j \theta ) \ . $$
Note also that to have a LPTI orbital symmetry we have to choose $\theta = \theta (x)$, and in this case $\la = (f \cdot \nabla) \theta$. \EOP

\subsection{Invariant trajectories}
\label{sec:IT}

We can ask the same question asked above for full symmetries also in the framework of orbital symmetries; that is, we can wonder which solutions to a given dynamical system have trajectories which are invariant under a given vector field. In this case, as we are in the framework of dynamical systems, the vector fields will be of the form
\beq \^v \ = \ \vphi^i (x,t) \, \frac{\pa}{\pa x^i} \ + \ \tau (x,t) \, \frac{\pa}{\pa t} \ ; \eeq
note that this is considerably more general than \eqref{eq:sVF} or \eqref{eq:sVFgen}, i.e. of those we actually want to consider here.

\medskip\noindent
{\bf Lemma 7.} {\it A trajectory $\ga$ of the dynamical system $\xd = f(x)$ is invariant under $\^v$ if and only if there is a smooth function $\theta (x,t)$ such that
\beql{eq:coll} \[ \vphi^i (x,t) \ - \ \theta (x,t) \ f^i (x) \]_\ga \ = \ 0 \ . \eeq}

\medskip\noindent
{\bf Proof.} A trajectory $\ga$ of the dynamical system is invariant under $\^v$ if and only if the latter and the dynamical vector field $\^f = f^i \pa_i$ are collinear on $\ga$; this is just the condition expressed in eq. \eqref{eq:coll}. \EOP

\medskip\noindent
{\bf Remark 5.} This result should be compared with Lemma 6 above. In fact, it says that on the invariant trajectory $\ga$ orbital symmetries reduce to the trivial ones identified in Lemma 6. \EOR

\section{Symmetry reduction for ODEs and symmetry adapted coordinates}

Let us briefly recall how the knowledge of $\G_f$, or part thereof\footnote{This is a relevant aspect of Lie theory: even if we are not able to determine the full symmetry algebra $\G_f$, we can anyway use any symmetry we have determined. See e.g. \gcite{\symref} for details.}, can be used to simplify the dynamical system \eqref{eq:DS}. If we know a single vector field $\^s = s(x) \nabla \in \^\G_f$, we can change coordinates and pass to new coordinates $(w^1,...,w^n)$ (note that our change of coordinates involve only the ``spatial'' ones) such that in these the symmetry vector field is written, say, as
\beql{eq:31} \^s \ = \ \frac{\pa}{\pa w^n} \ . \eeq That is, $\^s$ is now along one of the coordinate axes; the new coordinates are therefore \emph{symmetry adapted} ones.

\subsection{Symmetry reduction}
\label{sec:s3}

The commutation relation $[\^f,\^s] = 0$ being satisfied independently of the choice of coordinates, it will hold also in the new one; and in view of \eqref{eq:31}, this means that $\^f$ is independent of $w^n$. Thus in the new coordinates the system will actually read
\beql{eq:32} \dot{w} \ = \ g(w) \ , \ \ \ \ \pa g^i / \pa w^n \ = \ 0 \ \ \ \forall i=1,...,n \ . \eeq
In other words, we now have a system in (the first) $n-1$ coordinates, plus a last equation which tells how $w^n (t)$ evolves in time depending on the solution $\( w^1 (t) , ... , w^{n-1} (t) \)$ of that system. We will refer to this last equation as the \emph{reconstruction equation}; note that it actually amounts to an integration, i.e.
\beq w^n (t) \ = \ w^n (t_0) \ + \ \int_{t_0}^t g \[ w^1 (\theta ) , ... , w^{n-1} (\theta) \] \ d \theta \ . \eeq

\medskip\noindent
{\bf Example 1.} Consider the system
\beql{eq:33} \begin{cases} \xd \ = \ \a (r^2) \, x \ - \ \b (r^2 ) \, y &, \\
\yd \ = \ \b (r^2) \, x \ + \ \a (r^2 ) \, y &; \end{cases} \eeq
here $r^2 = x^2 + y^2$, and $\a,\b$ are smooth functions. This is symmetric under rotations, generated by the vector field $\^s = y \pa_x - x \pa_y$. Passing to polar coordinates $(r,\vartheta )$ we have $\^s = \pa_\vartheta$, and the time evolution is given by
\beql{eq:34} \begin{cases} \dot{r} \ = \ 2 \, r^2 \ \a (r^2) &, \\
\dot{\vartheta} \ = \ \b (r^2 ) &. \end{cases} \eeq
Rotation invariance is now expressed by r.h.s. terms not depending on $\vartheta$. \EOE

\medskip\noindent
{\bf Remark 6.} As shown by this Example, passing to symmetry adapted coordinates may involve a singular change of coordinates. \EOR
\bigskip

In the case we know a subalgebra $\G_0 \sse \G_f$, the procedure can be iterated $d$ times, where $d$ is the dimension of a maximal solvable subalgebra\footnote{Note that most often other symmetries, not belonging to $\G_0$,  will be lost in the procedure.} of $\G_0$ \gcite{Olv1}; the vector fields in $\G_0$ to be successively rectified should be chosen ``in the proper order'', i.e. the one dictated by the derived series of $\G_0$  \gcite{Olv1}. In this way \eqref{eq:DS} is reduced to a $(n-d)$ dimensional system plus $d$ ``reconstruction equations''.

If we know a constant of motion $P(x)$, we can also take advantage of this by passing to coordinates $z(x)$ such that say $z^n = P(x)$. In these new coordinates, eq. \eqref{eq:DS} will read
\beql{eq:35} \begin{cases} \dot{z}^i \ = \ h^i (z) & (i = 1,...,n-1) \ , \\
\dot{z}^n \ = \ 0 \ . & \end{cases} \eeq
Note that the $h^i$ do in general also depend on $z^n$.

\medskip\noindent
{\bf Example 2.} If in \eqref{eq:33} we choose $\a = 0$, then $r$ is a constant of motion. In this case the polar equations \eqref{eq:34} are just $\dot{\vartheta} = h (r^2 )$, $\dot{r} = 0$. \EOE
\bigskip

More in general, if we know $p$ functionally independent CM $\{ P_1 (x) , ... , P_p (x) \}$ we can pass to coordinates $\( z^1,...,z^{n-p}; \pi^1, ... , \pi^p \)$ with $\pi^i (x) = \pi^1 [ P_1 (x) , ... , P_p (x) ]$ (the simplest choice being of course $\pi^i = P_i (x)$); in these coordinates we have
\beql{eq:37} \begin{cases} \dot{z}^i \ = \ h^i (z,\pi) &, \\ \dot{\pi}^j \ = \ 0 &. \end{cases} \eeq

Note that the $\pi^j$ enter in \eqref{eq:37} only as parameters; the problem is therefore reduced to the study of a $(n-p)$ dimensional dynamical system with $p$ parameters. Fixing the value of the $\pi$'s identifies a submanifold $M_\pi \ss M$, by definition invariant under the flow of \eqref{eq:DS}; one is therefore legitimate to consider the reduction of \eqref{eq:DS} to this, i.e.
\beql{eq:38} \xd \ = \ f_{(\pi)} (x) \ , \ \ \ x \in M_\pi \ , \ \ \ f_\pi = f |_{M_\pi} \ : \ M_\pi \to \T M_\pi \ . \eeq

In the following, we will assume that the reduction corresponding to global CM and/or symmetries has already been performed.

\medskip\noindent
{\bf Remark 7.} The CM can be seen as invariant functions for the one-parameter group $\exp[ \la \^f]$. Then, under suitable assumptions of the action of this group, e.g. that it is compact, or acting properly\footnote{If we consider a linear vector field $\^f$ we have a subgroup of GL (n) and these conditions are satisfied.},  a celebrated theorem of Hilbert (for polynomial functions \gcite{Hil,OlvHil}; extended by Schwarz \gcite{Schw1,Schw2} to the smooth case) ensures \gcite{Mic1,Mic2,Mic3,Mic4,AbS,Sar} that there exists a basis $\{ \Theta_1 , ... , \Theta_d \}$ of invariant functions (respectively, polynomials) such that any $P(x) \in I_f$ (polynomial $P(x) \in I_f$) can be written as a function (as a polynomial) of the $\Theta$'s, i.e. as $P(x) \ = \ \wt{P} [ \Theta_1 (x) , ... , \Theta_d (x) ]$. \EOR

\medskip\noindent
{\bf Remark 8.} Note that it is not necessary that $\mathrm{dim} (M_\pi ) = [ \mathrm{dim} (M) - p ]$; e.g. in Example 2 above, for $r=0$ we have a zero-dimensional $M_\pi = \{ 0 \}$. \EOR

\medskip\noindent
{\bf Remark 9.} It is clear from \eqref{eq:37} that a symmetry vector field $\^s$ will be, in the $(z,\pi)$ coordinates, of the form
\beql{eq:39} \^s \ = \ h^i (z,\pi) \, \frac{\pa}{\pa z^i} \ + \ \eta^j (\pi) \, \frac{\pa}{\pa \pi^j} \ . \eeq
This can be checked by an explicit computation. \EOR

\medskip\noindent
{\bf Remark 10.} It may be worth remarking explicitly that $(f \cdot \nabla) P = 0$ together with $\{f,s\} = 0$ do \emph{not} imply $(s \cdot \nabla) P = 0$; see Example 3 below. \EOR

\medskip\noindent
{\bf Example 3.} Consider again \eqref{eq:33} with $\a (r^2) = 0$ and $\b (r^2) = 1$. In this case $r$ is a CM; the system is symmetric not only under rotations but also under the vector field $\^s = x \pa_x + y \pa_y$, which generates scaling transformations and do not leave $r$ invariant. \EOE

\subsection{Topology of trajectories and Lie-point symmetries}
\label{sec:s4}

It should be briefly recalled that LP transformations, being smooth and locally invertible\footnote{The invertibility is global if we consider proper -- that is, not just local -- groups of LP transformations.}, cannot transform a set $A \ss M$ into a topologically different one.

This applies in particular to the trajectory of a point $x \in M$ under the flow of \eqref{eq:DS}; so trajectories of ``topologically special'' types, e.g. fixed points, periodic orbits, or quasi-periodic ones filling densely a (topological) $k$-dimensional torus, cannot be transformed into orbits of different types. That is, isolated fixed points, isolated periodic orbits, isolated $k$-dimensional invariant tori are invariant under LP transformations \gcite{Cds,G93}.

\medskip\noindent
{\bf Remark 11.} Note that if were considering \emph{full solutions} rather than just their \emph{trajectories} -- or, in other words, the graph $\{ t , x(t)\}$ rather than the set $\{ x(t) , \ t \in R_+ \}$ -- there would be no topological difference between graphs of solutions corresponding to topologically different trajectories. In this sense, orbital symmetries turn out to be more relevant than full ones. \EOR

\subsection{Smooth structures and Lie-point symmetries}
\label{sec:s5}

It should be understood that preservation of topology of trajectories is not the only limitation on the type of solutions which can be connected by a LP transformation. Roughly speaking, preservation of topology just follows from the fact that LP transformations are one to one (invertible) and continuous; the fact they are not only $\mathcal{C}^0$ but $\mathcal{C}^\infty$ (or analytic), and hence not only continuous but also \emph{uniformly continuous}, poses further constraints \gcite{Cds,G93}, as we are going to discuss in this Section.

A first obvious remark is that a polynomial or $\mathcal{C}^\infty$ transformation can not map a $\mathcal{C}^\infty$ solution into a $\mathcal{C}^k$ one ($k$ finite); this will play a role in the discussion of center manifolds (see Section \ref{sec:s11} below).

Another useful remark is that a polynomial LP transformation cannot connect two trajectories which diverge exponentially. More in general, let $x(t)$ be the solution to \eqref{eq:DS} with initial datum $x(0) = x_0$; this identifies an invariant curve under $\^f$, i.e. the trajectory $\ga_0 \ss M$ of $x(t)$.

Suppose now that the linearization $A$ of $\^f$ at $x_0$ has positive eigenvalues in some direction $\xi \in \T_{x_0} M$, transversal to $\ga_0$; let $\^s$ be a vector field such that $\^s (x_0) = \xi$, and let $$ x_1  \ = \ e^{\eps \^s} \, x_0 $$ with $\eps$ small; denote by $\ga_1$ the trajectory of \eqref{eq:DS} issued from $x_1$. We consider an interval $\ga_0 (\ell)$ on $\ga_0$ of length $\ell > 0$, and a tubular neighborhood $u_\de$ of $\ga_0$, of radius $\de > 0$.

The $\ga_1$ will leave $u_\de$ for $\ell$ long enough, as the trajectories $\ga_0$ and $\ga_1$ diverge with a positive Lyapounov exponent \gcite{GuckHolm}. By the uniform continuity of LP transformations, $\ga_1$ and $\ga_0$ cannot be connected by a LP transformation.

For given $\de_0$ and $\ell$, by taking $\eps < \de_0$ small enough, we can guarantee that $|e^{\eps \^s} x - x| < \de_0$ $\forall x \in \ga_0 (\ell)$; but for any $\de_0$ and $\ell$ we can find a point $x^1 = e^{\eps \^s} x_0$, with $\eps < \de_0$, such that $|\ga_1 (\la) - \ga_0 (\la) | > \de_0$ for $\la_0 < \la < \ell$.

\medskip\noindent
{\bf Remark 12.} There is a point, in the above reasoning, that should be emphasized: the separation of the trajectories is exponential in the curvilinear coordinate along $\ga_0$\footnote{By taking $\de$ small enough, we can always choose a local system of coordinates in $u_\de$ such that $\ga_0$ corresponds to, say, $x^2=x^3=...=x^n = 0$, so the exponential in $x^1$ is well defined. (Notice for this argument to extend globally it is required that the group action is regular.)}, but this does not necessarily imply that two solutions running on these trajectories separate exponentially in time.

On the other hand, exponential separation in time of two solutions does not forbid that the corresponding trajectories are connected by a LP transformation, as Example 4 below shows. Also, the trajectories can separate exponentially in the curvilinear coordinate even if the solutions do not separate exponentially in time, see Example 5. \EOR

\medskip\noindent
{\bf Example 4.} Consider the simple system
$$ \begin{cases} \xd \ = \ x &, \\ \yd \ = \ y &; \end{cases} $$ this has solutions $x(t) = x_0 e^t$, $y(t) = y_0 e^t$, and trajectories are just straight half-lines $y=cx$ through the origin, covered at increasing speed. Consider two such motions, with initial data respectively $p_0 = (a,0)$ and $p_1 = (a,\eps)$. The LP transformation $e^{\eta \^s}$ with $\^s$ the vector field $\^s = [x \pa_y - y \pa_x] $ and $\eta = \arctan (\eps/a)$ transforms the trajectory $y=0$ issued from $p_0$ into the trajectory $y = (\eps/a) x$ issued from $p_1$. \EOE

\medskip\noindent
{\bf Example 5.} Consider the system
$$ \begin{cases} \xd \ = \ f(x) & \\
\yd \ = \ f(x) \cdot y &, \end{cases} $$
for which $dy/dx = y$, i.e. the trajectories are given by
$$ y \ = \ c \ e^x $$
with $c$ a constant (which can be zero). The trajectories $(x,0)$ issued from $p_0 = (a,0)$ and $(x,\eps e^x)$ issued from $p_1 = (a,\eps e^a)$ diverge exponentially in $t$. If $f(x)$ is such that $\exp[ x(t)]$ is not an exponential in $t$ (e.g. if $f(x) =\exp[ - x]$), the solutions do not diverge exponentially in time. \EOE

\section{Conditional symmetries. General setting}
\label{sec:s6}

We will now briefly recall the setting for the determination and use of \emph{conditional symmetries}, following \gcite{LeWin} (see \gcite{LWjmp,LRT2,LWB,LYW,PR1,PR2,PS1,PS2,CK,Olvp,Vor1,Vor2} for further detail and related topics). The special features appearing in the dynamical systems case will be discussed in the next Sect.\ref{sec:s7}, while in this Section, including Example 6 below, we will for a moment consider PDEs -- which were the original framework for conditional symmetries -- in order to recall the general theory.

Let us consider a differential equation $\De$ and its symmetry algebra $\G_\De$; we will denote by $\yb$ the independent variables and by $\ub$ the dependent ones. A generic vector field $\^\eta$ will be written in terms of these variables as
\beql{eq:61} \^\eta \ = \ \xi^i (\yb,\ub) \, \frac{\pa}{\pa y^i} \ + \ \vphi^a (\yb,\ub) \, \frac{\pa}{\pa u^a} \ . \eeq
If $\De$ is of order $n$, the vector fields $\^\eta \in \G_\De$ are those satisfying
\beql{eq:62} \begin{cases} \^\eta^{(n)} \cdot \De \ = \ 0 & \\
\De \ = \ 0 & , \end{cases} \eeq where $\^\eta^{(n)}$ is the $n$-th prolongation of $\^\eta$ \gcite{\symref}.

As already recalled, knowledge of $\G_\De$ allows for a reduction of the equation $\De$; in particular, one can look for solutions $\ub = f (\yb)$ which are invariant under a subgroup $\G_0 \sse \G_\De$. With this invariance \emph{ansatz}, eqs.\eqref{eq:62} reduce to simpler ones. Indeed, now we can express $\ub (\yb)$ in terms of the differential invariants of $\G_0$ \gcite{Olv1}.

It should be noted that given a vector field $\^\eta_0$ of the form
\beql{eq:61b} \^\eta_0 \ = \ \xi_0^i (\yb,\ub) \, \frac{\pa}{\pa y^i} \ + \ \vphi_0^a (\yb,\ub) \, \frac{\pa}{\pa u^a} \ , \eeq
the condition of invariance of $\ub = f (\yb)$ under $\^\eta$ reads
\beql{eq:63} \De_0 \ := \ \vphi_0^a \ - \ \xi_0^i \ \( \frac{\pa u^a}{\pa y^i} \) \ = \ 0 \ , \eeq
so that \emph{if} we know \emph{apriori} that $\^\eta_0 \in \G_\De$, and therefore that \eqref{eq:62} is satisfied, the solutions to $\De$ which are invariant under $\^\eta_0$ can be seen as the solutions of the system made by \eqref{eq:62} and \eqref{eq:63}, i.e. of the system
\beql{eq:64} \begin{cases} \De \ = \ 0 &, \\ \De_0 \ = \ 0 &. \end{cases} \eeq

The simple key observation here is that it may also happen that a solution to $\De$ is invariant under a vector field which is \emph{not} in $\G_\De$. This suggests that symmetry reduction could be possible, and useful, also considering symmetries (of solutions) which are not in $\G_\De$ \gcite{LeWin}.

\medskip\noindent
{\bf Definition 3.} {\it The vector field $\^\eta_0$ is a \emph{conditional symmetry} for the equation $\De$ if and only if there are solutions to this which are invariant under $\^\eta_0$.}

\medskip\noindent
{\bf Example 6.} Consider the PDE for $u = u(x,y) \in R$ written as
\beql{eq:65} u_{xx} \ + \ u_{yy} \ = \ - \, 2  \, u \ + \ \a \, (x^2 \ + \ \a \, y^2) \, u \ - \  ( 1 \, - \, \a ) \, y \, u_y \ , \eeq
where $\a$ is a real constant, $\a \not= 1$. This is \emph{not} invariant for rotations in the $(x,y)$ plane, generated by the vector field $\^\eta_0 = y \pa_x - x \pa_y$. On the other hand, the function
\beql{eq:66} f(x,y) \ = \ \exp \[ - \frac12 \, \( x^2 + y^2 \) \] \eeq
is rotationally symmetric, and $u = f(x,y)$ is a solution to \eqref{eq:65}. \EOE
\bigskip

We stress that one is by no means guaranteed that there exist nontrivial solutions to $\De$ which are invariant under a given vector field which is in $\G_\De$; in any case, the theory of conditional symmetries aims at detecting solutions which are invariant under vector fields which are not in $\G_\De$, so that one often considers vector fields $\^\eta \in \G_\De$ as trivial conditional symmetries.

Solving the original equation $\De = 0$ with the \emph{ansatz} of invariance under $\^\eta_0$ amounts to solving \eqref{eq:64}. Note that there we do not care about the first equation in the system \eqref{eq:62}. If we give a generic $\^\eta_0$, anyway, the system \eqref{eq:64} will in general have no solution, so that this method is useful \emph{only} if we are able to determine the $\De_0$ compatible with the original $\De$, i.e. the $\^\eta_0$ which can yield symmetric solutions to $\De$.

A method for determining these, and the corresponding invariant solutions, does indeed exist \gcite{Ovs,LYW,PatWin}, and we now briefly illustrate it. We stress that albeit we always (for ease of notation) discuss invariance under a single vector field, the whole discussion is immediately extended to a Lie algebra of vector fields, simply by considering invariance under all of its generators at the same time.

Let us now consider $\^\eta_0$ mentioned above to be not specified; that is, the functions $\xi^i_0 (\yb,\ub)$ and $\vphi^a_0 (\yb,\ub)$ are undetermined. The solutions to $\De$ which are invariant under $\^\eta_0$ will still be given by solutions to \eqref{eq:64}; however this should now be seen not as a system of equations for just the unknown function $f(\yb)$ determining $\ub = f (\yb)$, but as a system of equations for the unknown functions $\xi^i_0 (\yb,\ub)$, $\vphi^a_0 (\yb,\ub)$, and $f (\yb)$.

We can apply to \eqref{eq:64} the known methods for solving a system of PDEs (we note that even in the case where $\De$ is an ODE,  $\De_0$ is a first order PDE); in particular, we can apply symmetry methods.

Let us first of all determine the symmetries of \eqref{eq:64}; the determining equations for this can be written as
\beql{eq:67} \begin{cases}
\^\eta^{(n)} \cdot \De \ = \ 0 & \\
\^\eta^{(1)} \cdot \De_0 \ = \ 0 & \\
\De \ = \ 0 & \\
\De_0 \ = \ 0 &. \end{cases} \eeq

Notice that if we choose to consider $\^\eta_0 = \^\eta$ which is a full symmetry for the equation $\Delta$ (which is surely possible, as $\^\eta_0$ is completely generic) the second of these is automatically satisfied, as $\De_0 = 0$. This means that \eqref{eq:67} is equivalent to \eqref{eq:62} and \eqref{eq:63}, i.e. it yields no restriction on $\^\eta_0$.

We can now apply to \eqref{eq:64} the usual symmetry reduction method, and obtain (symmetric) solutions in this way. Such solutions to \eqref{eq:64} (if they exist) give us a vector field $\^\eta_0$ and at the same time a function $f$, i.e. a solution $\ub = f (\yb)$ invariant under $\^\eta_0$; by construction, $\ub = f (\yb)$ is also a solution to the original equation $\De = 0$. Again by construction, $\^\eta_0$ is a symmetry of \eqref{eq:64}, and if it exists (which we assume from now on, not to repeat over and over this specification) it is thus obtained by solving the corresponding determining equations \eqref{eq:67}.

We can solve the determining equations \eqref{eq:67} by the well known algorithms for solving determining equations \gcite{Olv1,Olv2,Ovs,Ste} (also by computer algebra if needed \gcite{symCA}); once this is done, we can choose a specific $\^\eta_0$ among the solutions to \eqref{eq:67}, and pass to consider the corresponding eq.\eqref{eq:64}, i.e. determine the $\^\eta_0$-invariant solutions to the original equation.

The vector fields $\^\eta_0$ which solve \eqref{eq:67} are symmetries of the system \eqref{eq:64} but in general (except for $\^\eta_0 \in \G_\De$, which represent here the trivial set) \emph{not} of $\De = 0$ alone. They are therefore called \emph{conditional symmetries} of $\De$, as they are symmetries of $\De$ when this is subject to the additional condition \eqref{eq:63} \gcite{LeWin}.

\medskip\noindent
{\bf Remark 13.}
Note that the additional (or side) condition \eqref{eq:63}, i.e. $\De_0$, depends on $\^\eta_0$ itself. Therefore the conditional symmetries of a given equation $\De$ do \emph{not} in general form an algebra, as they are ordinary symmetries of \emph{different} systems. Note also that the set of conditional symmetries of $\De$ does naturally carry an action of $\G_\De$, hence  conditional symmetries for a given equation can be subdivided into conjugacy classes for this action, thus leading to a classification.  \EOR
\bigskip

It should be stressed that albeit the introduction of conditional symmetries was motivated by the search for invariant solutions, we can very well have vector fields which leave invariant \emph{no} solutions to $\De$, but which transform a subset of solutions into the same subset. Such vector fields, to be formally defined in a moment, are then said to be \emph{partial symmetries} for $\De$. We will not discuss them in detail here, albeit later on we will briefly refer to the possibility of extending our discussion to such symmetries; for details about partial symmetries the reader is referred to \gcite{CGpart,Cpart}.

\medskip\noindent
{\bf Definition 4.} {\it The vector field $\wt{\eta}_0$ is a \emph{partial symmetry} for the equation $\De$ if and
only if there is a subset $\mathcal{S}_0$ of solutions to $\De$ which is globally invariant under $\wt{\eta}_0$.}

\medskip\noindent
{\bf Remark 14.} Proper symmetries and conditional symmetries are extreme (degenerate) cases of partial symmetries: for proper symmetries the subset $\mathcal{S}_0$ coincides with the set of \emph{all} solutions to $\De$, while for conditional symmetries the set $\mathcal{S}_0$ reduces to a single solution (or union of solutions, each of them)  individually invariant under $\wt{\eta}_0$). \EOR

\medskip\noindent
{\bf Remark 15.} Together with the ``direct problem'', i.e. determining the possible conditional symmetries of a given equation $\De$ (that is, the possible symmetries of its solutions beside the symmetries of the equation itself), it is of obvious physical interest also the ``inverse problem'': that is, given a symmetry group $G$ determining which equations can admit solutions with the $G$ symmetry. This problem has been tackled by Levi, Rodriguez and Thomova \gcite{LRT} and by Pucci and Saccomandi \gcite{PScond1,PScond2}. We will not discuss it here. \EOR

\medskip\noindent
{\bf Remark 16.} The relations of the Levi-Winternitz theory of conditional symmetries with the Michel theory on the geometry of group action are briefly commented upon in \gcite{GLW}. \EOR

\section{Conditional symmetries and conditional constants \\ of motion for dynamical systems}
\label{sec:s7}

In this Section we specialize the general theory discussed in the previous Section to the case of dynamical systems. We will also discuss the relation between existence of conditional (LTPI) symmetries and existence of ``conditional'' constants of motion.

\subsection{Conditional symmetries for dynamical systems}

Let us now see how the general discussion of conditional symmetries given in the previous Section specializes in the case of Dynamical Systems. In this case, $\De$ will be a system of first order autonomous ODEs as in \eqref{eq:DS},
\beql{eq:71} \De \ \equiv \ \ \xd^i \ - \ f^i (\xb) \ = \ 0 \ ; \eeq
here, as usual, $x \in M \sse R^N$, $f : M \to \T M$. Correspondingly, we write generic LP vector fields on $M \times R$ as
\beql{eq:72} \^\eta \ = \ \s^i (\xb,t) \, \pa_i \ + \ \tau (\xb , t) \, \pa_t \eeq
with of course $\pa_i = \pa / \pa x^i$.

The equation $\De_0$, see eq.\eqref{eq:63}, is in this case
\beql{eq:73} \s^i (\xb,t) \ - \ \tau (\xb,t) \, \xd^i \ = \ 0 \ ; \eeq
using \eqref{eq:71} this reads simply
\beql{eq:74} \s^i (\xb,t) \ - \ \tau (\xb,t) \ f^i (\xb) \ = \ 0 \ . \eeq

We want again to focus on autonomous vector fields which do not act on $t$ as well, i.e. look for Lie-Point Time-Independent (LPTI) symmetries, or actually LPTI conditional symmetries; these will be called for short \emph{configurational symmetries}\footnote{This denomination has a double advantage: on the one hand it stresses that the vector fields depends only on the configuration of the system (and not on the time at which it is reached); on the other hand it recalls the concept of configurational invariants \gcite{SLC} which we will meet soon.}. Thus we consider vector fields of the form
\beql{eq:76} \^\eta \ = \ s^i (\xb) \, \pa_i \ . \eeq

Now the equation $\^\eta^{(1)} \cdot \De = 0$ gives, before restriction to solution of $\De$,
\beql{eq:77} \xd^j \, \pa_j s^i \ - \ s^j \, \pa_j f^i \ = \ 0 \ , \eeq
while a solution $x(t) = \xi (t)$ is invariant under \eqref{eq:76} if
\beql{eq:78} \De_0 \ := \ \ s^i [\xi_1 (t),...,\xi_n (t)] \ = \ 0 \ . \eeq
This states that the LPTI vector fields must vanish altogether on the solution.

We could now repeat the discussion of Sect.\ref{sec:s6}; this would just require to consider \eqref{eq:71}, \eqref{eq:77} and \eqref{eq:78} rather than \eqref{eq:61} and \eqref{eq:63}. We will not bore the reader with such a repetition, and just introduce some definitions which are natural in view of it.

\medskip\noindent
{\bf Definition 5.} {\it The vector field \eqref{eq:76} is a (LPTI) \emph{conditional symmetry} for the dynamical system \eqref{eq:71} if and only if there are solutions $x^i = \xi^i (t)$ to this, such that \eqref{eq:78} is satisfied for all $t$.}

\subsection{Conditional orbital symmetries for dynamical systems}
\label{sec:COSDS}

We recall now that (see Sect.\ref{sec:OS}) the vector field $\^\eta$ is an \emph{orbital symmetry} \gcite{WalOS,WalMPS} of a dynamical system if it maps solutions trajectories into solution trajectories. We may extend the concept of conditional symmetries to this kind of symmetries.

\medskip\noindent
{\bf Definition 6.} {\it The vector field $\^\eta_0$ is a \emph{conditional orbital symmetry} for a dynamical system \eqref{eq:DS} if there are solution trajectories of the dynamical system which are invariant under $\^\eta_0$.}

\medskip\noindent
{\bf Remark 17.} If we do not ask for invariance  under \eqref{eq:76} of the full solution $x = \xi (t)$, but only of its \emph{trajectory}
\beql{eq:79} \ga_\xi \ := \ \{ \xi (t) \ , \ t \in R \} \ , \eeq
the invariance equation is \eqref{eq:coll} discussed above (see Lemma 6 and Lemma 7). In the LPTI case this reads simply
\beql{eq:710} \[ \vphi^i (x) \ - \ \theta (x) \ f^i (x) \]_\ga \ = \ 0 \ . \eeq
This requires, indeed, that the (orbital) symmetry vector fields and the dynamical one are collinear on the considered trajectory $\gamma$. \EOR

\medskip\noindent
{\bf Remark 18.} A conditional symmetry leaves some solution invariant; this implies the trajectory is also invariant, and hence a conditional symmetry is also a (special case of) conditional orbital symmetry. \EOR

\medskip\noindent
{\bf Remark 19.} Note that one may have cases in which no single trajectory is invariant under $\^\eta_0$, but a subset\footnote{We restrict to a subset because if the whole set of trajectories is globally invariant, we have a proper orbital symmetry.} of them is left globally invariant; in this case we will speak of \emph{partial orbital symmetries}, due to similarity with the partial symmetries mentioned above \gcite{CGpart,Cpart}.

An example of this case can be built as follows. Let $\{ 0 \}$ and $r_1$ be the origin and the unit circle in $R^2$. Consider the  dynamical system in $R^2 \backslash \{ 0 \} \backslash r_1$ given in polar coordinates $(\varrho,\vartheta)$ by
$$ \dot\varrho \ = \ \varrho \, (1 - \varrho) \ + \ K (\varrho) \ \a (\varrho,\vartheta) \ , \ \ \ \dot\vartheta \ = \ \om \ + \ K(\varrho) \b (\varrho,\vartheta)  \ ,  $$ where $K(\varrho)$ is a smooth function exactly vanishing for $\varrho > 1 $ (but not for $\varrho < 1$) and $\a,\b$ are smooth but sufficiently ``strange'' functions. Then the inward spiralling trajectories living in $\varrho > 1$ are not rotationally invariant but are mapped one into the other by rotations, while the solutions living in $\varrho < 1$ are not rotationally invariant and in general not mapped one into the other. \EOR
\bigskip

We now note that the discussion of the previous Section \ref{sec:s6} shows that (LPTI) conditional symmetries of $\De$ are also obtained as standard (LPTI) symmetries of the system \beql{eq:711} \begin{cases} \De \ = \ 0 &, \\ \De_0 \ = \ 0 &, \end{cases} \eeq where $\De_0$ is given by the invariance condition \eqref{eq:63}. It is easily seen that the same holds for conditional orbital symmetries, provided $\De_0$ is now given instead by \eqref{eq:coll}; or in the LPTI case by \eqref{eq:710}. This suggest an alternative definition of conditional orbital symmetries:

\medskip\noindent
{\bf Definition 6'.} {\it The vector field \eqref{eq:76} is a \emph{conditional orbital symmetry} for the dynamical system \eqref{eq:71} if and only if it is a standard orbital symmetry of the system \eqref{eq:711}.}

\medskip\noindent
{\bf Remark 20.} It should be stressed that Definition 6 and Definition 6' are equivalent, while of course the notions of conditional symmetry and conditional orbital symmetries are different (with the former implying the latter, see Remark 18), as the foregoing Example 7 does clearly show. \EOR

\medskip\noindent
{\bf Example 7.} Consider the dynamical system
\beql{eq:exCSDS} \begin{cases} \xd \ = \ \a (x,y) \, x \ - \ \b (x,y) \, y &,\\ \yd \ = \ \b (x,y) \, x \ + \ \a (x,y) \, y &. \end{cases} \eeq
(Note that here $\a,\b$ are not necessarily functions of $r^2 = x^2+y^2$, at difference with Example 1 above.) To this is associated the vector field
\beq \^f \ = \ \( \a (x,y) \, x \ - \ \b (x,y) \, y \) \, \pa_x \ + \ \( \b (x,y) \, x \ + \ \a (x,y) \, y \) \, \pa_y \ . \eeq  Consider the rotation vector field
\beq \^\eta_0 \ = \ - y \, \pa_x \ + \ x \, \pa_y \ . \eeq

\medskip\noindent
{\bf (a)} Choose $\a (x,y) = 1 - x^2 - y^2$, $\b (x,y) = \om \not= 0$; the trajectories (apart from that starting and remaining in the origin) are spiralling towards the unit circle $r^2 = 1$ both from outside and from inside. None of these is rotationally invariant but those starting on the unit circle itself; obviously a rotation maps a trajectory into a trajectory in all cases. As for full solutions, again rotations generated by $\^\eta_0$ do not leave these invariant (apart from the trivial one) invariant, but map any solution into a different solution.

\medskip\noindent
{\bf (b)} Consider general functions $\a,\b$ \emph{not} being a function of $r^2$ alone.
It results
\begin{eqnarray*}
\[ \^\eta_0 , \^f \] &=& \( x^2 \, \a_y \ - \ x \, y \, \( \a_x \, + \, \b_y \) \ + \ y^2 \, \b_y \) \, \pa_x \\
&=& \ + \ \( x^2 \, \b_y \ + \ x \, y \, \( \a_y \, - \, \b_x \) \ - \ y^2 \, \a_y \) \, \pa_y \ . \end{eqnarray*}
Thus, as expected, $\^\eta_0$ is not a symmetry unless
$$ \a (x,y) \ = \ a (r^2 ) \ , \ \ \b (x,y) \ = \ b (r^2 ) \ , \ \ \ r^2 = x^2 + y^2 \ . $$

Consider now, with $\rho$ the unit circle, the special case where $\a$ and $\b$ satisfy
\beql{eq:exCSDSab} \[ \a (x,y) \]_\rho \ = \ 0 \ , \ \ \ \[ \b (x,y) \]_\rho \ = \ \om \ . \eeq
In this case the circle $\rho$ is an invariant manifold; the restriction of our dynamical system to $\rho$ is given by
\beq \begin{cases} \xd \ = \ - \ \om \, y &,\\ \yd \ = \ \om \, x  &, \end{cases} \eeq
i.e. by uniform rotations with angular speed $\om$. Note again that none of these solutions is invariant under $\^\eta_0$, which just maps them one into the other; thus (if we decide to disregard the trivial solution) $\^\eta_0$ is not a conditional symmetry for the dynamical system \eqref{eq:exCSDS} when \eqref{eq:exCSDSab} are satisfied; on the other hand, it is a conditional orbital symmetry. \EOE

\subsection{Conditional constants of motion for dynamical systems}
\label{sec:CCMDS}

In the same spirit, it is natural to consider, besides conditional symmetries and conditional orbital symmetries, also \emph{conditional constants of motion}; these will be defined more precisely in a moment. As recalled in the Introduction they were introduced by Sarlet, Leach and Cantrijn \gcite{SLC} under the name of \emph{configurational invariants}; they are also implicitly considered in \gcite{LeWin} (see eq. (2.9) therein), provided one restricts the general discussion by Levi and Winternitz to the case of dynamical systems instead on considering the general PDE case. For applications of conditional constants of motion, see e.g. \gcite{PucRos1,PucRos2,PucRos3}.

\medskip\noindent
{\bf Definition 7.} {\it A function $\mu : \wt{M} \to R$, where $\wt{M} \sse M$, is a conditional constant of motion for the dynamical system \eqref{eq:DS} if and only if there are level sets of $\mu$ which are invariant under the dynamical vector field $\^f = f^i (x) \pa_i$.}

\medskip\noindent
{\bf Remark 21.} We stress that we do not require $\mu$ to be defined in all of $M$; thus e.g. if the level set $\mu^{-1} (0)$ is the union of several disjoint sets, we can restrict consideration to $\wt{M}$ a tubular neighborhood of one of these disjoint sets, and the invariance of one of the disjoint component of $\mu^{-1} (0)$ is enough for $\mu$ to qualify as a conditional symmetry. \EOR

\medskip\noindent
{\bf Remark 22.} Any full constant of motion is also a conditional constant of motion; we will see these as trivial conditional constants of motion, and focus on ``proper'' conditional constants of motion (omitting to specify this at each step). \EOR

\medskip\noindent
{\bf Remark 23.} One could provide equivalent definitions in the case of $\mu : \overline{M} \to R$; the details of such an extension are left to the reader. \EOR
\bigskip

We would like to characterize conditional constants of motion in the same way as we did for conditional symmetries and conditional orbital symmetries. Indeed, in the same way as a conditional symmetry is a symmetry of the system \eqref{eq:711}, we have that

\medskip\noindent
{\bf Definition 7'.} {\it A function $P(\xb )$, $P: M \to R$, is a \emph{conditional constant of motion} for $\De$ if and only if it is an ordinary constant of motion for the system \eqref{eq:711}.}
\bigskip

Finally, it should be noted that by definition \eqref{eq:711} admits solutions only for $\De_0$ corresponding to a conditional symmetry; thus Definition 7 makes sense only in this case. That is, \emph{conditional constants of motion are associated to conditional symmetries} (and to Cicogna's partial symmetries \gcite{CGpart,Cpart}) \emph{and actually to conditional orbital symmetries}, as we are going to discuss in greater detail in the following.

Before doing this, it is worth considering again -- now from the present point of view -- simple examples given in previous Sections.

\medskip\noindent
{\bf Example 8.} Consider again the simple systems of Examples 1--3, say in $M = R^2 \backslash \{ 0 \}$ to discard trivial solutions, namely
\beql{eq:712} \begin{cases} \xd \ = \ \a (r^2) \, x \ - \ \b (r^2) \, y &, \\
\yd \ = \ \b (r^2) \, x \ + \ \a (r^2) \, y &; \end{cases} \eeq
in polar coordinates this is simply $\dot{r} = 2 r^2 \a(r)$, $\dot\vartheta = \b (r^2)$.

If $\a (r^2)$ has some nontrivial zero, say $r = r_0$, then the solutions on the invariant circle $r_0$, i.e. the solutions to the system \eqref{eq:711} where $\De$ is given by \eqref{eq:712} and $\De_0$ by $r=r_0$, are exchanged among themselves by rotations, i.e. by actions of the vector field $\^\eta = y \pa_x - x \pa_y = \pa_\vartheta$.

Note that while this shows that $r^2 = x^2+y^2$ is a conditional constant of motion in the sense of Definition 7, and that $\pa_\vartheta$ is a conditional orbital symmetry in the sense of Definition 6, we cannot conclude that $\pa_\vartheta$ is a conditional symmetry in the sense of Definition 5.  In fact, as already mentioned, rotations do exchange solutions on this circle among themselves so that $\pa_\vartheta$ is a \emph{partial symmetry} of $\De$ in the sense of Cicogna \gcite{CGpart,Cpart} but not a conditional symmetry; on the other hand, the vector field
\beql{eq:714} \^\eta \ := \ \b (r_0^2) \, \pa_t \ + \ \[ y \, \pa_x \ - \ x \, \pa_y \] \ = \ \b (r_0^2) \, \pa_t \ + \ \pa_\vartheta \eeq
is a conditional symmetry for \eqref{eq:712} if $r_0$ is such that $\a (r_0^2) = 0$. \EOE

\medskip\noindent
{\bf Example 9.} Consider $M = R^3$ with coordinates $(x,y,z)$, and write $r^2 = x^2 + y^2$. Consider the dynamical system
\beql{eq:715} \begin{cases} \xd \ = \ \a(r^2) \, x \ - \ \b (r^2) \, y & \\
\yd \ = \ \b (r^2) \, x \ + \ \a (r^2) \, y & \\
\zd \ = \ f(z) \ + \ g(x,y) \, r^2 & \end{cases} \eeq with $\a,\b,f,g$ smooth functions of their arguments.

It is immediate to check that $\pa_\vartheta = y \pa_x - x \pa_y$ is now a conditional symmetry, corresponding to solutions on the $z$ axis, i.e. such that $x(t)=0=y(t)$. Note also that, unless $g(x,y) = \wt{g} (x^2+y^2)$, we have that $\pa_\vartheta$ is not an ordinary symmetry of \eqref{eq:715}; note also that it is an ordinary symmetry for \eqref{eq:711}. \EOE
\bigskip

\medskip\noindent
{\bf Remark 24.} The previous Example 8 shows that requiring the conditional symmetries of a dynamical system to be time-independent can be too restrictive a condition, especially in view of Definition 6; thus one should consider LP vector field of general form. However, conditional symmetries of a dynamical system \emph{can} also be time-independent, as shown by Example 9. \EOR

\bigskip

The discussion above -- see in particular Examples 8 and 9 -- clearly shows why in the sense of Definition 6 standard symmetries \emph{can} also be conditional symmetries, but are not necessarily so.

\section{Conditional symmetries, conditional constants of motion,  and invariant manifolds}
\label{sec:s8}

It is clear from the discussion of Sect.\ref{sec:s7} that a close connection exists between conditional symmetries and conditional constants of motion of a dynamical system on one hand, and its \emph{invariant manifolds} on the other hand. We are now going to discuss this relation in some detail.

Let us first consider LPTI conditional symmetries; in this case the additional condition $\De_0$ in \eqref{eq:711} is given by \eqref{eq:78}. In this way, any vector field identifies a manifold $S_0 \sse M$, and solutions to \eqref{eq:711} are nothing else than the solutions to $\De$ which lie entirely in $S_0$. These are associated to the set $\Xi_0 \sse S_0$ of points $\xb$ for which $\dot\xb \in \T S_0$; that is,
\beql{eq:81} \Xi_s \ = \ \{ x \in S_0 \ : \ \xd = f (x) \in \T_x S_0 \} \ . \eeq
It is clear that $\Xi_s$ is an invariant manifold for the dynamical system \eqref{eq:71}; that is (by definition)
\beql{eq:82} f \ : \ \Xi_s \ \to \ \T \Xi_s \ . \eeq

Note that if we have a different conditional symmetry with a vector field $\^\eta'$ identifying $S_0' \sse M$ and the same $\Xi_{s'} = \Xi_s$ (this is the case e.g. for $\^\eta' = k \^\eta$), this would yield the same invariant solutions as $\^\eta$. This fact, and ease of notation, suggest to consider in particular the case where $\Xi_s = S_0$, i.e. the case where $S_0$ is itself an invariant manifold for the dynamical system:
\beql{eq:83} f \ : \ S_0 \ \to \ \T S_0 \ . \eeq

In both cases \eqref{eq:82} and \eqref{eq:83}, the determination of conditional symmetries greatly simplifies -- or even solves -- the problem of determining invariant manifolds for the dynamical system under study. It should be stressed that this simple remark puts at once at our disposal the powerful and \emph{completely algorithmical} methods developed for the study of symmetries and conditional symmetries of differential equations, and allows to use them for the determination of invariant manifolds of dynamical systems.

Let us now concentrate on the case \eqref{eq:83}. The invariant manifolds being determined by \eqref{eq:78}, it follows at once that the $s^i (x)$ are conditional constants of motion for the dynamical system $\xd^i = f^(x)$ under study. More in general, any function $P(x)$ such that $S_0$ is a level set for $P$ will be a conditional constant of motion for the dynamical system. Similarly, in case \eqref{eq:82} any function $P(x)$ such that $\Xi_s$ is a level set for $P$ will be a conditional constant of motion for the dynamical system. In this case we can first restrict $P(x)$ to $S_0$, call $P_0 (x)$ the restricted function, and then consider simply the level sets of $P_0$.

\medskip\noindent
{\bf Remark 25.} Clearly, if $P(x)$ is a conditional constant of motion for a given dynamical system, so is also $P_c (x) := P(x) + c$, where $c$ is a constant. Denoting by $I_0$ the algebra of constant functions on $M$, this suggests that we should consider CM and conditional constants of motion modulo $I_0$. \EOR
\bigskip

Note now that $S_0$ is a level set for $P(x)$ if and only if
\beql{eq:84} \( {\bf v} (x) \cdot \nabla \) \, P(x) \ = \ 0 \eeq
for \emph{any} vector field $\^v = v^i (x) \pa_i$ such that $\^v : S_0 \to \T S_0$; the latter condition just means that
\beql{eq:85} v^j(x) \pa_j s^i (x) \ = \ 0 \ \ \ \forall i \ \ \ \mathrm{on} \ \ \ {\bf s} (x) \ = \ 0 \ . \eeq
We have therefore a simple method to determine the $P(x)$ associated to a given conditional symmetry $\^\eta$. Indeed, once $s^i (x)$ are known, eq.\eqref{eq:85} is a linear equation for the $v^j$,
\beql{eq:86} A^i_{\ j} (x) \ v^j (x) \ = \ 0 \eeq
where the matrix $A$ is defined by
\beql{eq:87} A^i_{\ j} (x) \ = \ \( \frac{\pa s^i}{\pa x^j} \)_{S_0} \ . \eeq

Once the general solution ${\bf v} (x)$ of this is known, we can pass to determine $P(x)$ using \eqref{eq:84}, i.e. the characteristic equation
\beql{eq:88} \frac{d x^1}{v^1 (x)} \ = \ ... \ = \ \frac{d x^N}{v^N (x)} \ . \eeq

A close connection also exists between conditional orbital symmetries and conditional constants of motion. Indeed, if $\^s$ is a conditional orbital symmetry, the (union of) trajectories invariant under it provide invariant manifolds; such trajectories correspond, once $s(x)$ is fixed, to the manifold $\chi_s$ of points satisfying \eqref{eq:79}, i.e.
\beql{eq:89} \chi_s \ := \ \{ x \in M \ : \ \ |f(x)| \cdot |s(x) | \ = \ \( f(x) , s(x) \) \} \ . \eeq

Note that $S_0 \sse \chi_s$; this corresponds to the fact that for dynamical systems, any conditional symmetry is a conditional orbital symmetry (see Remark 18). We can then proceed as in \eqref{eq:84}--\eqref{eq:88}, with the role of $S_0$ played by $\chi_s$.

\subsection{Partial symmetries and conditional constants of motion}

As already mentioned, it is also interesting to consider \emph{partial symmetries} \gcite{CGpart,Cpart}, and correspondingly \emph{partial orbital symmetries}; these have a relation with conditional constants of motion.

Indeed, let $\^\eta = s^i (x) \pa_i$ be a partial symmetry, and let us consider the corresponding set $\S_0$ (see Definition 7) of invariant solutions; as $f$ and $s$ are smooth functions, the trajectories of these span a manifold which we also call $S_0$. By construction this is invariant under the flow of $f$ (and of $\^\eta$, of course); we can then apply again the discussion between \eqref{eq:84} and \eqref{eq:88}.

In this way a partial symmetry is also characterized as a vector field $\^\eta$ such that it exists a manifold $S_0 \sse M$ for which, denoting by $\rho_0$ the operator of restriction to $S_0$,
\begin{eqnarray}
& & \rho_0 \[ \{ f , s \} \] \ = \ 0 \ , \nonumber \\
& & \rho_0 (\^f ) : S_0 \to \T S_0 \ , \ \ \rho_0 (\^\eta ) : S_0 \to \T S_0 \ . \label{eq:810} \end{eqnarray}

Finally, let us consider a partial orbital symmetry, as defined in Sect.\ref{sec:s7}; again the \emph{closure} of the invariant set of trajectories $\mathcal{T}_0$ (i.e. the union of points $x \in \ga$ for $\ga \in \mathcal{T}_0$) defines a manifold in $M$ which is invariant for both $\^f$ and $\^\eta$, and the same considerations apply.

\subsection{Module structure of conditional symmetries}
\label{sec:s9}

The relation between conditional symmetries and conditional constants of motion and between orbital symmetries and conditional constants of motion discussed in the previous Section \ref{sec:s8} also means that there is a module structure of conditional symmetries over the corresponding conditional constants of motion, pretty much analogous to the structure of $\G_\De$ as a module over $I_\De$ discussed in Section \ref{sec:s2}, as we now show.

Let us reverse the point of view of Section \ref{sec:s8}, and consider manifolds $S_0 \sse M$ identified by
\beql{eq:91} s^i (x) \ = \ 0 \ \ \ \forall i=1,...,n \ . \eeq

Suppose we have chosen $s^i (x)$ such that $S_0$ is invariant under $\^f$; then any vector field $\^\eta : S_0 \to \T S_0$ is a partial symmetry (possibly a conditional symmetry) for $\De$. Moreover $\^s = s^i \pa_i$ is then a conditional symmetry for $\De$. More in general, we could have several vector fields which are conditional symmetries for $\De$ and which leave invariant one (or more, or all) solution lying in $S_0$.

It should be stressed that, once $S_0$ has been fixed, the associated partial symmetries form a Lie algebra, which we will denote as $\G_{S_0}^W$; this follows form the definition of partial symmetries  (see Definition 4), or equivalently from \eqref{eq:810}.

As for conditional symmetries, those which leave invariant \emph{some} solutions lying in $S_0$ (not the same solution for different conditional symmetries) do \emph{not} form an algebra. On the other hand, if we fix a \emph{given} set of solutions lying in $S_0$ (possibly all of them), then the conditional symmetries leaving these invariant \emph{do} form an algebra. To be specific, let us consider the conditional symmetries leaving invariant all the solutions in $S_0$, and let us call their algebra $\G_{S_0}^C$; clearly,
\beql{eq:92} \^s \in \G^C_{S_0} \sse \G_{S_0}^W \ . \eeq

Let us also consider the set $I_S$ of conditional constants of motion associated to $S_0$, i.e. the smooth functions $P : M \to R$ for which $S_0$ is a level set (the relation between $I_{S_0}$ and $S_0$ is discussed in Sect.\ref{sec:s8}); these do form an algebra. It is then immediate to remark, as in Sect.\ref{sec:s2}, that $\G_{S_0}^C$ and $\G_{S_0}^W$ are not only algebras, but also the structure of a Lie module over $I_{S_0}$.

Similar remarks, and corresponding results, hold when considering configurational symmetries.

\section{Special points and manifolds}
\label{sec:specialpts}

In this final Section, we will discuss the relations between the different objects we have been considering -- that is, conditional symmetries, orbital symmetries, conditional orbital symmetries, and conditional constants of motion -- and the canonical invariant manifolds associated to fixed points studied in Dynamical Systems theory \gcite{HPS,Ruelle,GuckHolm,Ver1,ArnEMS}.

\medskip\noindent
{\bf Remark 26.} We will only develop consequences of our discussion above, and not aim to discuss in general the relation between symmetries and invariant manifolds for dynamical systems (which would require a review in itself). This is a classical topic, see e.g. the early works by Steeb, Wulfman and Bluman  \gcite{Steeb1,Steeb2,Wul,Blu} as well as more recent contributions \gcite{CG94a,FGG,CGWjlt} (also, closing the circle, inspired by Olver and Rosenau's approach to invariant solutions for PDEs through ``side conditions'' \gcite{OlvRos1,OlvRos2}).

\subsection{Hyperbolic fixed points, stable and unstable manifolds}
\label{sec:s10}

In the last two sections we have discussed the relation between conditional symmetries and conditional constants of motion on the one hand, and invariant manifolds on the other hand; this was considered under different points of view, which are natural for the investigation of symmetry properties. We want now to reverse our point of view, and consider these relations, and the involved symmetry properties altogether, from the point of view of dynamical systems theory. This discussion is strongly related to the one developed in the context of ordinary symmetries in \gcite{CG94a}, to which we refer for further details in that context.

A number of results are already well known in this direction, but these  mainly concern ordinary rather than conditional symmetries \gcite{Cds,CGihp,CG93,CG94a,CG94b,CG94c}; we briefly recall them, referring to the original papers for details and proofs.

First of all, as already remarked in Sect.\ref{sec:s4}, if the smooth dynamical system \eqref{eq:DS} admits an \emph{isolated} fixed point $\xb_0$,
\beql{eq:102} f^i (\xb_0 ) \ = \ 0 \ \ \ \ i=1,...,n \ , \eeq
then any LPTI (ordinary) symmetry $\^\eta = s^i (x) \pa_i$ must satisfy
\beql{eq:103} s^i (\xb_0) \ = \ 0 \ \ \ \ i = 1,...,n \ . \eeq
Note that \eqref{eq:103} is not a necessary condition for $\^\eta$ to be a conditional symmetry; on the other hand, it is a sufficient one -- the invariant solution being simply $\xb (t) = \xb_0$.

Let us now suppose $\xb_0$ is moreover an \emph{hyperbolic} fixed point \gcite{HPS,Ruelle,GuckHolm,Ver1,ArnEMS} (this also implies it is an isolated one); this means that all the eigenvalues of the linear operator
\beql{eq:104} A \ := \ (D f) (\xb_0) \eeq
have nonzero real part. In this case one associates to $\xb_0$ the unique local stable and unstable manifolds $W_s$ and $W_u$; these are invariant manifolds under the flow of $\^f$ (see e.g. \gcite{HPS,Ruelle,GuckHolm,Ver1,ArnEMS} for the role they play in dynamical systems theory).

The situation described above for the fixed point case is, with obvious adaptations, met again here. That is, it can be proved \gcite{CGihp,CG93,CG94a} that a necessary (but not sufficient) condition for $\^\eta = s^i (x) \pa_i$ to be an ordinary symmetry of \eqref{eq:DS} is that
\beql{eq:105} \^\eta \ : \ W_s \to \T W_s \ , \ \ \ \^\eta \ : \ W_u \to \T W_u \ . \eeq
This implies that the restrictions of $\^\eta$ to $W_s$ and $W_u$ are well defined, and indeed it is immediate to see that a condition (again, necessary but not sufficient) for $\^\eta$ to be an ordinary symmetry is that
\beql{eq:106} \[ \^\eta , \^f \]_{W_s} \ = \ 0 \ , \ \ \ \[ \^\eta , \^f \]_{W_u} \ = \ 0 \ . \eeq

Conversely, each of this is clearly a sufficient -- but not necessary -- condition for $\^\eta$ to be a conditional symmetry\footnote{At first sight, one may think this just guarantees $\^\eta$ is a partial symmetry, but \eqref{eq:106} also implies \eqref{eq:103}.} for \eqref{eq:DS}.

Moreover, the invariant manifolds $W_s$, $W_u$ will be identified by equations
\beql{eq:107} w_s (\xb ) \ = \ 0 \ , \ \ \ w_u (\xb ) \ = \ 0  \eeq for certain smooth function $w_s$ and $w_u$. Note that these are not unique, albeit $W_s$ and $W_u$ are: any function with the same zero level set would do. We also note that functions vanishing on $W_s$ (respectively, on $W_u$) form an algebra, which we denote as $I_W^s$ (respectively, $I_W^u$).

It is clear that the $w_s (\xb )$, $w_u (\xb )$ are conditional constants of motion for \eqref{eq:DS}. As discussed above, the functions satisfying either condition in \eqref{eq:106} will form a Lie module over the algebra $I_W^s$ and respectively $I_W^u$.

These considerations extend \emph{a fortiori} to conditional configuration symmetries.

It should also be mentioned that the property of being invariant under any LPTI (ordinary) symmetry is not peculiar to stable and unstable manifolds, but extends to any transversally hyperbolic manifold \gcite{CG94a}.

\subsection{Non hyperbolic fixed points, center manifolds}
\label{sec:s11}

The case where the fixed point $\xb_0$ is not hyperbolic is slightly more complicated. In this case (see \gcite{CG94a} for precise conditions and statements of results), together with stable and unstable manifolds one has to consider center manifolds $W_c$ \gcite{HPS,Ruelle,GuckHolm,Ver1,ArnEMS}.

Contrary to what happens for stable and unstable manifolds, the center manifold in \emph{not} uniquely defined; moreover even in the case the dynamical system is $C^\infty$, the center manifold could be non-smooth, or only smooth of class $C^k$ \gcite{ArnEMS}. Anyway, each of the center manifolds is still an invariant manifold for $\^f$; nevertheless, the results holding in the hyperbolic case do not extend immediately to this setting.

The obstacle to such an extension is intimately related to the non-uniqueness of the center manifold, i.e. to terms beyond all orders in the perturbation expansion of the center manifold (these are in turn related to resurgent functions \gcite{ArnEMS,Ecalle}).

We will therefore consider only the \emph{Poincar\'e-Dulac center manifold} $W_c^0$; this is defined to be the (analytic) center manifold constructed perturbatively or, equivalently, the (infinite order) jet of center manifolds \gcite{ArnODE,ArnEMS,ArnGM}. Thus we are essentially disregarding (setting to zero) all the terms beyond all orders in perturbation: all the center manifolds differ for non-perturbative terms and thus share the same jet at all orders.

\medskip\noindent
{\bf Example 10.} We give a simple example to illustrate the problem with center manifolds. Consider the dynamical system
\beql{eq:DS1} \begin{cases} \xd \ = \ - \, x^3 &, \\
\yd \ = \ - \, y &. \end{cases} \eeq
It is immediate to see that all the curves
\beql{eq:DS2} y \ = \ \a \ e^{- 1/x^2} \eeq
are center manifolds, for any $\a \in R$. The perturbation expansion for all these is anyway the same, and just yields $y = 0$.

Note that albeit all these center manifolds are $C^\infty$, the one corresponding to $\a = 0$ is the only one to be analytic. \EOE
\bigskip

If we consider the Poincar\'e-Dulac center manifold, the same results given in the previous Sect.\ref{sec:s10} apply \gcite{CG93,CG94a}, in particular for what concerns eqs. \eqref{eq:105}, \eqref{eq:106} and \eqref{eq:107}, up to a simple rephrasing amounting to consider $W_c^0$ instead of $W_s$ or $W_u$.

We stress that in general, i.e. for generic $W_c \not= W_c^0$, even $\^\eta : W_c \to \T W_c$ does \emph{not} hold, and we can only be sure that a LPTI symmetry transforms a solution lying on a center manifold $W_c^a$ into a solution lying on a center manifold $W_c^b$, where $W_c^a$ and $W_c^b$ may happen to coincide but in general are different; see again \gcite{CG93,CG94a}.

\subsection{Bifurcation of fixed points}
\label{sec:s12}

In the case where the dynamical system \eqref{eq:DS} depends on a real control parameter $\la$, so that we actually deal with
\beql{eq:121} \xd^i \ = \ f^i (\xb; \la ) \ , \eeq
and $\xb_0$ is a fixed point for all values of $\la$ (at least within an interval of interest),
\beql{eq:122} f^i (\xb_0 ; \la ) \ = \ 0 \ \ \ \forall \la \eeq
which is undergoing a simple bifurcation \gcite{Ruelle,GuckHolm,ChH,GSS} at $\la = \la_0$, the setting and results of the previous Sections \ref{sec:s10} and \ref{sec:s11} combine nicely.

Let us suppose, to fix ideas, that $\xb_0$ is stable for $\la < \la_0$ and only a real eigenvalue, or a pair of complex conjugate ones, crosses the imaginary axis with positive speed at $\la = \la_0$. Then, beside stable manifolds $W_s (\la)$, the system has a center manifold $W_c$ for $\la = \la_0$, and an unstable manifold $W_u (\la)$ for $\la > \la_0$. The theorem on persistence of transversally hyperbolic manifolds \gcite{Ruelle} ensures that the center manifold $W_c$ can be uniquely defined by the requirement to be the limit of $W_u (\la)$ for $\la \to \la_0$ (from the right).

Thanks to this property in this case, again, we do no meet the problems related to non-uniqueness of the center manifold; moreover, as $W_u (\la)$ corresponds to its perturbation expansion, the limit $W_c$ is precisely the Poincar\'e-Dulac center manifold. Once again we refer to \gcite{CG94a} for details.

\section{Discussion}
\label{sec:disc}

In this Section we will present a short discussion, in the form of some remarks, in particular concerning limitations (and hence possible extensions) of our discussion. 

\medskip\noindent
{\bf Remark 27.} As customary in the symmetry theory of differential equations (including dynamical systems) we made no difference between conservative -- in particular, Hamiltonian -- and dissipative dynamical systems. In fact, we have always worked in the general case, and in general considered systems for which Energy is not conserved; moreover, even in the case where Energy is conserved, this is -- from the point of view of our discussion -- on equal footing with pother conserved quantities (if any).
 
It goes without saying that in the Hamiltonian case one can perform a Legendre transformation and pass to a Lagrangian description; in this case symmetry considerations can to a large extent be encompassed in the formalism of the (general) Noether theorem \gcite{Olv1,Koss}. We believe the present treatment is already too long without entering into the details of the special features of this case, and have thus confined ourselves to the general case, referring the reader to the works mentioned above \gcite{Olv1,Koss} for the conservative framework. In this respect, see also \gcite{MSS,MarmoB,CIMM,Sarda,Urban}. \EOR

\medskip\noindent
{\bf Remark 28.} Needless to say, in between the cases of autonomous dynamical systems (which we have been considering here) and general time-dependent ones, one could consider the special case of time-periodic ones. In general terms, it would be natural in this case to restrict consideration to Lie-point  symmetries which are themselves time-periodic (rather than time-independent).  

From the point of view of symmetries, periodic systems are characterized by a \emph{discrete} symmetry, corresponding to time translation by the period $T$. This is \emph{not} a Lie-point symmetry, hence this case would not add any special feature in our general discussion. On the other hand, when we consider conditional symmetries (and conditional constants of motion), these are -- or at least can be -- related to discrete symmetries. We hope some of the readers of the present discussion can investigate periodic systems from this point of view. 

One may note, in this context, that from the point of view of conditional constant of motions (and possibly conditional full or partial integrability) reflection symmetries appear to be more fruitful. This can be extended to reversing \gcite{Sevryuk,Birkhoff,Bibikov,Lamb1,Lamb2,RobQui} and $k$-reversing dynamical systems as well \gcite{Lamb3,Lamb4}. Applications of these to the search of special solutions of dynamical systems is discussed e.g. in \gcite{QuiRob,LambNic,BLH}; see also e.g. \gcite{BPV,Dias} for applications in Fluid Dynamics. \EOR

\medskip\noindent
{\bf Remark 29.} Coming back to time-periodicity, one should recall that symmetries, in particular discrete symmetries, may enforce the appearance of time-periodic solutions in autonomous systems (this is also related to the spontaneous breaking of continuous symmetries). In this respect see, among others, \gcite{RodVan,Mont,CraKno,Swift,Robbins,Zhil,Ruck,VerGal}. \EOR

\medskip\noindent
{\bf Remark 30.}  Orbital symmetries \gcite{WalOS,WalMPS,Wal19} have played a relevant role in our discussion. One should note that they have already been used (since their introduction) to pinpoint special solutions to dynamical systems. One can see in this regard -- and also about other features of orbital symmetries -- e.g. \gcite{WalOS,WalMPS,Wal19,CGWjlt,FGG,KinWal,SchroWal,GSW,SchWal,HadWal,CLPW,Goeke}. The latter references deal with chemical reactions equations; these are also discussed, from a point of view relevant to our discussion, in \gcite{Broad1,Broad2} (see also \gcite{Broad3}). \EOR

\medskip\noindent
{\bf Remark 31.} We have been considering ``full'' dynamical systems, without explicit mention of perturbation theory and perturbative approach, except for some mention in Section \ref{sec:s11}. When dealing with dynamical systems admitting a stationary solution (an \emph{equilibrium}), it is entirely natural to consider (Poincar\'e) \emph{normal forms} around this. (This also extends to systems admitting periodic or multi-periodic solutions and normal forms around these.) This calls immediately for the analysis of symmetries, orbital symmetries, conditional symmetries etc. for systems which are in normal form. This topic has of course been considered in the literature, and actually it plays a prominent role in several of the references provided in our discussion above in this Section. We refrain from entering even a cursory discussion of this topic and the related literature. \EOR 

\medskip\noindent
{\bf Remark 32.} The reader has probably remarked that we only presented very simple Examples; their goal was indeed to illustrate in the simplest possible terms the concepts and results discussed in the paper. On the other hand, the (rather ample) bibliography given in the course of our discussion contains a wealth of applications in many different fields; we preferred (also to avoid making this too long paper even longer) giving indications of the original papers rather than discussing some examples taken from these. \EOR

\section{Conclusions}

We have considered the well known notion (due to Levi and Winternitz) of  \emph{conditional symmetries}. Our discussion also naturally called for consideration of an extension of conditional symmetries, i.e. that (due to Cicogna) of \emph{partial symmetries}.

We focused on the applications of these notions, originally developed for determining special solutions to PDEs, to ODEs; and focused in particular on their applications to Dynamical Systems.

In particular, we have discussed the relations of conditional and partial symmetries with two other established notions (these confined by their nature to Dynamical Systems), i.e. \emph{orbital symmetries} and \emph{conditional constants of motion} (also known as \emph{configurational invariants}, discussing the interrelations among these.

We also looked in more detail, in Section \ref{sec:specialpts}, at how the main characters of our discussion enter in a classical topic, i.e. the relations between symmetry properties and invariant sets (points, trajectories, manifolds) for dynamical systems.

Albeit we provided no new results, we trust that our discussion clarifies on the one hand the interrelations between these notions, and on the other hand how concepts -- such as that of conditional symmetries -- created in the analysis of PDEs can also be of interest in the seemingly simpler arena of analysis of Dynamical Systems.

In Section \ref{sec:disc} we have also discussed a number of related topics we have \emph{not} dealt with, including relevant applications, and provided some references for these.


\addcontentsline{toc}{section}{\ \ \ \ Appendix A.  Simple examples of orbital symmetries}

\section*{Appendix A.  Simple examples of orbital symmetries}

Orbital symmetries have played a relevant role in our discussion. As remarked in Section \ref{sec:OS}, these are less widely known than standards symmetries, so we provide here some extremely simple examples to help the reader who has not met these before fix ideas. More substantial examples are to be found in the literature, see in particular  \gcite{WalOS,WalMPS,Wal19,CGWjlt,FGG,KinWal,SchroWal,GSW,SchWal,HadWal,CLPW,Goeke}.

\medskip\noindent
{\bf Example A1.} Consider first the system 
\beql{eq:OSexa1} \begin{cases} \xd \ = \ - \om (x^2 + y^2) \ y &,\\ \yd \ = \ \om (x^2 + y^2) \ x &,\end{cases} \eeq
with $\om$ an arbitrary smooth function satisfying $\om (0) = 0$; 
its solutions are given by uniform circular motions preserving the radius $r = \sqrt{x^2 + y^2}$; circles are travelled on with constant angular speed $\om (r^2)$. Note this system corresponds to the vector field
\beql{eq:OSexa1VF} X_1 \ = \ \om (x^2+y^2) \ \[ - \, y \ \pa_x \ + \ x \, \pa_y \] \ . \eeq
The rotation vector field
\beql{eq:OSexaXr} X_r \ = \ - y \, \pa_x \ + \ x \, \pa_y \eeq maps solutions into solutions and is a proper symmetry (and also an orbital one), while the scaling vector field
\beql{eq:OSexaXs} X_s \ = \ x \, \pa_x \ + \ y \, \pa_y \eeq does \emph{not} map solutions into solutions (unless $\a$ reduces to a constant), but it maps solution trajectories (circles) into solution trajectories (circles, again). \EOE

\medskip\noindent
{\bf Example A2.} 
Consider now the system 
\beql{eq:OSexa2} \begin{cases} \xd \ = \ - \b (x,y) \ y &,\\ \yd \ = \ \b (x,y) \ x &,\end{cases} \eeq with $\b$ an arbitrary smooth function satisfying $\b (0,0) = 0$. 
Its solutions are given by motions on circles, thus preserving the radius $r = \sqrt{x^2 + y^2}$; circles are now in general travelled on with varying angular speed $\b (x,y)$, i.e. with an angular speed which depends in an arbitrary (smooth) way on both $r$ and $\theta$. Note this system corresponds to the vector field
\beql{eq:OSexa2VF} X_2 \ = \ \b (x,y) \ \[ - \, y \ \pa_x \ + \ x \, \pa_y \] \ . \eeq In this case both $X_r$ and $X_s$ defined above are in general not proper symmetries, and they are both orbital symmetries. \EOE

\medskip\noindent
{\bf Example A3.} For the examples \eqref{eq:OSexa1} and \eqref{eq:OSexa2} considered above, it is straightforward to compute commutators of relevant vector fields, i.e. \eqref{eq:OSexa1VF}, \eqref{eq:OSexa1VF}, with the orbital symmetry ones, i.e. \eqref{eq:OSexaXr} and \eqref{eq:OSexaXs}. 

We obtain (omitting functional dependencies for ease of notation)
\begin{eqnarray*} 
\[ X_r , X_1 \] \ = \  0 & ; & \ \[ X_s , X_1 \] \ = \ \( \frac{2 \, (r^2) \, \om' }{\om} \) \ X_1 \ ; \\
\[ X_r , X_2 \] \ = \ \( \frac{x \, \b_y \ - \ y \, \b_x}{\b} \) \ X_2 & ; & \ \[ X_s , X_2 \] \ = \ \( \frac{x \, \b_x \ + \ y \, \b_y}{\b} \) \ X_2 \ . \end{eqnarray*}
These confirm the result of Lemma 4 above (see Section \ref{sec:OS}). 
We also note, in view of Lemma 5 (see again Section \ref{sec:OS}), that $[X_r,X_s] = 0$. \EOE

\addcontentsline{toc}{section}{\ \ \ \ Acknowledgements}

\subsection*{Acknowledgements}

I thank Sebastian Walcher (RWTH Aachen) for a critical reading of a preliminary version of this work and for suggesting several improvements \& corrections. My work is partially supported by the project {\it ``Mathematical Methods in Non-Linear Physics''} (MMNLP) of INFN (CNS4), and by GNFM-INdAM. Most of this work was performed while enjoying the warm hospitality of SMRI.
The final version was improved thanks to the constructive criticism by an unknown Referee.

\label{lastpage}


\newpage

\addcontentsline{toc}{section}{\ \ \ \ References}

\begin{thebibliography}{129}

\bibitem{LeWin} D. Levi and P. Winternitz, ``Non-classical symmetry reduction: example of the Boussinesq equation'', {\it J. Phys. A} {\bf 22} (1989), 2915-2924

\bibitem{LWR} D. Levi and P. Winternitz. ``Continuous symmetries of difference equations'', {\it J. Phys. A} {\bf 39} (2005), R1-R63

\bibitem{WalOS} S. Walcher, ``Orbital symmetries of first order ODEs'', in {\it Symmetry and perturbation theory -- SPT98} (A. Degasperis and G. Gaeta eds.), World Scientific 1999, pp. 96--113

\bibitem{WalMPS} S. Walcher, ``Multi-parameter symmetries of first order ordinary differential equations'', {\it J. Lie theory} {\bf 9} (1999), 249-269

\bibitem{SLC} W. Sarlet, P.G.L. Leach and F. Cantrijn, ``First integrals versus configurational invariants and a weak form of complete integrability'', {\it Physica D} {\bf 17} (1985), 87-98



\bibitem{Olv1} P.J. Olver, {\it Application of Lie groups to differential equations}, Springer 1986

\bibitem{Olv2} P.J. Olver, {\it Equivalence, Invariants, and Symmetry}, Cambridge UP  1995

\bibitem{Ovs} L.V. Ovsijannikov, {\it Group analysis of differential equations}, Academic Press 1982; based on {\it Group properties of differential equations}, Novosibirsk 1962

\bibitem{Ste} H. Stephani, {\it Differential equations. Their solution using symmetries}, Cambridge University Press 1989

\bibitem{AVL} D.V. Alexseevsky, A.M. Vinogradov and V.V. Lychagin, {\it Basic Ideas and Concepts of Differential Geometry}, Springer 1991

\bibitem{CGb} G. Cicogna and G. Gaeta, {\it Symmetry and perturbation theory in nonlinear dynamics}, Springer 1999

\bibitem{KrV} I.S. Krasil'schik and A.M. Vinogradov, {\it Symmetries and conservation laws for differential equations of mathematical physics}, A.M.S. 1999

\bibitem{Win1} P. Winternitz, ``What is new in the study of differential equations by symmetry methods'', Group Theoretical Methods in Physics (Proceedings XV ICGTMP), R. Gilmore ed., World Scientific (1987).

\bibitem{Win2} P. Winternitz, ``Lie groups and solutions of nonlinear partial differential equations''; in: A. Ibort and M.A. Rodriguez eds., {\it Integrable systems, quantum groups, and quantum field theories}, Springer 1993, pp 429-495



\bibitem{ArnODE} V.I. Arnold, {\it Ordinary Differential Equations}, Springer 1992

\bibitem{ArnEMS} V.I. Arnold and Yu.S. Ilyashenko, ``Ordinary Differential Equations'', in \emph{Encyclopaedia of Mathematical Sciences, vol.1 (Dynamical Systems I)}, D.V. Anosov and V.I. Arnold eds., Springer 1998

\bibitem{GuckHolm} J. Guckenheimer and Ph. Holmes, {\it Nonlinear oscillations, dynamical systems, and bifurcation of vector fields}, Springer 1983

\bibitem{Ver1} F. Verhulst, {\it Nonlinear differential equations and dynamical systems}, Springer 1990

\bibitem{Cds} G. Cicogna, ``Lie-point symmetries and dynamical systems'', pp. 147-153, in  {\it Modern Group Analysis: Advanced Analytical and Computational Methods in Mathematical Physics} ( N. H. Ibragimov, M. Torrisi, A. Valenti eds.), Springer  1993

\bibitem{CGihp} G. Cicogna and G. Gaeta, ``Lie-point symmetries in bifurcation problems'',  Annales de l'IHP Physique th\'eorique. Vol. 56. No. 4. 1992.

\bibitem{CGncb} G. Cicogna and G. Gaeta, ``On Lie point symmetries in mechanics'', Il Nuovo Cimento B (1971-1996) 107 (1992): 1085-1096

\bibitem{G93} G. Gaeta, ``Autonomous systems, dynamical systems, LPTI symmetries, topology of trajectories, and periodic solutions'', {\it Int. J. Theor. Phys.} {\bf 32} (1993), 191-199

\bibitem{Kir} A.A. Kirillov, {\it Elements of the theory of representations}, Springer 1976

\bibitem{Hil} D. Hilbert, {\it Theory of algebraic invariants} (translated by B. Sturmfels and R.C. Laubenbacher), Cambridge UP 1994

\bibitem{OlvHil} P.J. Olver, {\it Classical Invariant Theory}, Cambridge UP 1999

\bibitem{Schw1} G. Schwartz, ``Smooth functions invariant under the action of a compact Lie group'', {\it Topology} {\bf 14} (1975), 63-68

\bibitem{Schw2} G. Schwartz, ``Lifting smooth homotopies of orbit spaces'', {\it Publ. Math. IHES} {\bf 51} (1980), 37-135


\bibitem{Mic1} L. Michel, ``Points critiques des fonctions invariantes sur une G-vari\'et\'e'',
{\it C. R. Acad. Sc. Paris} {\bf 272} (1971), 433-436

\bibitem{Mic2} L. Michel, ``Nonlinear group action. Smooth action of compact Lie groups on manifolds'', in {\it Statistical Mechanics and Field Theory}, R.N. Sen and C. Weil eds., Israel University Press 1971

\bibitem{Mic3} L. Michel, ``Symmetry defeects and broken symmetry. Configurations. Hidden symmetry'', {\it Rev. Mod. Phys.} {\bf 52} (1980), 617-653

\bibitem{Mic4} L. Michel and B.I. Zhilinskii, ``Symmetry, invariants, topology. Basic tools'', {\it Phys. Rep.} {\bf  341}  (2001), 11-84



\bibitem{AbS} M. Abud and G. Sartori, ``The geometry of spontaneous symmetry breaking'', {\it Ann. Phys.} {\bf 150} (1983), 307-372

\bibitem{Sar} G. Sartori, ``Geometric invariant theory: a model-independent approach to spontaneous symmetry and/or supersymmetry breaking'', {\it Rivista del Nuovo Cimento} {\bf 14.11} (1991), 1-120





\bibitem{LWjmp} D. Levi and P. Winternitz,
``Symmetries and conditional symmetries of differential-difference equations'',
{\it J. Math. Phys.} {\bf 34} (1993), 3713-3730

\bibitem{LRT2} D. Levi, M.A. Rodriguez, and Z. Thomova, ``The discretized Boussinesq equation and its conditional symmetry reduction'', {\it J Phys. A} {\bf 53} (2020), 045201

\bibitem{LWB} D. Levi, R. Rebelo and P. Winternitz eds., {\it Symmetries and Integrability of Difference Equations}, Springer 2017


\bibitem{LYW} D. Levi, R. Yamilov and P. Winternitz, {\it Continuous Symmetries and Integrability of Discrete Equations}, AMS 2023



\bibitem{PR1} P.J. Olver and Ph. Rosenau, ``The construction of special solutions to partial differential equations'', {\it Phys. Lett. A} {\bf 114} (1986), 107-112

\bibitem{PR2} P.J. Olver and Ph. Rosenau, ``Group-invariant solutions of differential equations'', {\it SIAM J. Appl. Math.} {\bf 47} (1987), 263-278

\bibitem{PS1} E. Pucci and G. Saccomandi, ``On the weak symmetry groups of partial differential equations'', {\it J. Math. Anal. Appl.} {\bf 163} (1992), 588-598

\bibitem{PS2} E. Pucci, ``Similarity reduction of partial differential equations'', {\it J. Phys. A} {\bf 25} (1992), 2631-2640

\bibitem{CK} P.A. Clarkson and M.D. Kruskal, ``New similarity reduction of the Boussinesq equation'', {\it J. Math. Phys.} {\bf 30} (1989), 2201-2213

\bibitem{Olvp} P.J. Olver, ``Direct reduction and differential constraints'', preprint 1993

\bibitem{Vor1} E.M. Vorob'ev, ``Partial symmetries of systems of differential equations'', {\it Soviet Math. Dokl.} {\bf 33} (1986), 408-412

\bibitem{Vor2} E.M. Vorob'ev, ``Reduction and quotient equations for differential equations with symmetries'', {\it Acta Appl. Math.} {\bf 23} (1991), 1-24

\bibitem{PatWin} J. Patera and P. Winternitz, ``Subalgebras of real three and four dimensional Lie algebras'', {\it J. Math. Phys.} {\bf 18} (1977), 1449-1455

\bibitem{symCA} L. Amata, F. Oliveri, and E. Sgroi, ``Optimal systems of Lie subalgebras: A computational approach'', {\it J. Geom. Phys.} {\bf 204} (2024), 105290

\bibitem{CGpart} G. Cicogna and G. Gaeta, ``Partial Lie-point symmetries of differential equations'', {\it J. Phys. A} {\bf 34} (2001), 491-512

\bibitem{Cpart} G. Cicogna, ``Partial symmetries and dynamical systems'', {\it Math. Meths. Appl. Sci.} {\bf 39}  (2016), 4171-4180

\bibitem{LRT} D. Levi, M.A. Rodriguez, and Z. Thomova, ``Differential equations invariant under conditional symmetries'', {\it J. Nonlin. Math. Phys.} {\bf 26} (2019), 281-293

\bibitem{PScond1} E. Pucci and G. Saccomandi, ``Partial differential equations admitting a given nonclassical point symmetry'', {\it Studies Appl. Math.} {\bf 145} (2020), 81-96

\bibitem{PScond2} E. Pucci and G. Saccomandi, ``Using symmetries \`a rebours'', {\it Studies Appl. Math.} {\bf 146} (2021), 99-117

\bibitem{GLW} G. Gaeta, ``On the conditional symmetries of Levi and Winternitz'', {\it J. Phys. A} {\bf 23} (1990), 3643-3645

\bibitem{PucRos1} K. Rosquist and G. Pucacco, ``Invariants at fixed and arbitrary energy. A unified geometric approach'', {\it J. Phys. A} {\bf 28} (1995), 3235-3252

\bibitem{PucRos2} G. Pucacco and K. Rosquist, ``Configurational invariants of Hamiltonian systems'', {\it J. Math. Phys.} {\bf 46} (2005), 052902

\bibitem{PucRos3} G. Pucacco and K. Rosquist, Energy dependent integrability, {\it J. Geom. Phys.} {\bf 115} (2017), 16-27

\bibitem{HPS} M. Hirsch, D. Pugh and M. Shub, {\it Invariant manifolds} (LNM 583), Springer 1975

\bibitem{Ruelle} D. Ruelle, {\it Elements of differentiable dynamics and bifurcation theory}, Academic Press 1989

\bibitem{Steeb1} W.H. Steeb, ``Non-linear dynamic systems, limit cycles, transformation groups, and perturbation techniques'', {\it J. Phys. A} {\bf 10} (1977), L221-L223

\bibitem{Steeb2} W.H. Steeb, ``Nonlinear autonomous dynamic systems, limit cycles, and one-parameter groups of transformations'', {\it Lett.  Math. Phys.} {\bf 2} (1977), 171-174

\bibitem{Wul} C.E. Wulfman, ``Limit cycles as invariant functions of Lie groups'', {\it J. Phys. A} {\bf 12} (1979), L73-L75

\bibitem{Blu} G. Bluman, ``Invariant solutions for ordinary differential equations'', {\it SIAM J. Appl. Math.} {\bf 50} (1990), 1706-1715

\bibitem{CG94a} G. Cicogna and G. Gaeta, ``Symmetry invariance and centre manifolds for dynamical systems'', {\it Nuovo Cimento B} {\bf 109} (1994), 59-76

\bibitem{FGG} E. Freire, A. Gasull, and A. Guillamon, ``Limit cycles and Lie symmetries'', {\it Bull. Sci. Math.} {\bf 131} (2007), 501-517

\bibitem{CGWjlt} G. Cicogna, G. Gaeta and S. Walcher, ``Side conditiona for ordinary differential equations'', {\it J. Lie Theory} {\bf 25} (2015), 125-146

\bibitem{OlvRos1} P.J. Olver and P. Rosenau, ``The construction of special solutions to partial differential equations'', {\it Phys. Lett. A} {\bf 114}  (1986), 107-112

\bibitem{OlvRos2} P.J. Olver and P. Rosenau, ``Group-invariant solutions of differential equations'', {\it SIAM J. Appl. Math.} {\bf 47} (1987), 263-278
    
\bibitem{CG93} G. Cicogna and G. Gaeta, ``Nonlinear Lie symmetries in bifurcation theory'', {\it Phys. Lett. A} {\bf 172} (1993), 361-364

\bibitem{CG94b} G. Cicogna and G. Gaeta, ``Poincar\'e normal forms and Lie point symmetries'', {\it J. Phys. A} {\bf 27} (1994), 461-476

\bibitem{CG94c} G. Cicogna and G. Gaeta, ``Normal forms and nonlinear symmetries, {\it J. Phys. A} {\bf 27} (1994), 7115-7124

\bibitem{Ecalle} J. Ecalle, ``Singularit\'es non abordables par la g\'eom\'etrie'', {\it Ann. Inst. Fourier} {\bf 42} (1992), 73-164

\bibitem{ArnGM} V.I. Arnold, {\it Geometrical methods in the theory of ordinary differential equations}, Springer 1983

\bibitem{ChH} S.N. Chow and J.K. Hale, {\it Elements of bifurcation theory}, Springer 19821


\bibitem{GSS} M. Golubitsky, D. Schaeffer and I. Stewart,   {\it Singularities and groups in bifurcation theory}, Springer 2012




\bibitem{Koss} Y. Kossmann-Schwarzbach, {\it Noether Theorems: Invariance and Conservation Laws in the 20th Century}, Springer, 2009

\bibitem{MSS} G. Marmo, E.J. Saletan and A. Simoni, ``A general setting for reduction of dynamical systems'', {\it J. Math. Phys.} {\bf 20} (1979), 856-860

\bibitem{MarmoB} G. Marmo, E.J. Saletan, A. Simoni and B. Vitale, {\it Dynamical systems. A differential geometric approach to symmetry and reduction}, Wiley 1985 

\bibitem{CIMM} J.F. Carinena, A. Ibort, G. Marmo and G. Morandi, {\it Geometry from dynamics, classical and quantum}, Springer 2015

\bibitem{Sarda} G. Sardanashvily, {\it Noether's Theorems, Applications in Mechanics and Field Theory}, Springer 2016

\bibitem{Urban} Z. Urban, F. Bajardi and S. Capozziello, ``The Noether-Bessel-Hagen symmetry approach for dynamical systems'', {\it Int. J. Geom. Meth. Mod. Phys.} {\bf 17} (2020), 2050215



\bibitem{Sevryuk} M.B. Sevryuk, {\it Reversible systems} (LNM Vol. 1211),  Springer 2006

\bibitem{Birkhoff} G.D. Birkhoff, ``Dynamical systems with two degrees of freedom'', {\it Proc. Nat. Acad. Sci.} {\bf 3} (1917), 314-316

\bibitem{Bibikov} Yu.N. Bibikov, ``On the stability of the zero solution of a periodic reversible second-order differential equation'', {\it Vestnik St. Petersb. Univ. Math.} {\bf 55} (2022), 297-300

\bibitem{Lamb1} J.S.W. Lamb, ``Reversing symmetries in dynamical systems'', {\it J. Phys. A} {\bf 25} (1992), 925-937

\bibitem{Lamb2} J.S.W. Lamb and J.A.G. Roberts, ``Time-reversal symmetry in dynamical systems: a survey'', {\it Physica D} {\bf 112} (1998), 1-39

\bibitem{RobQui} J.A.G. Roberts and G.R.W. Quispel, ``Chaos and time-reversal symmetry. Order and chaos in reversible dynamical systems'', {\it Phys. Rep.} {\bf 216} (1992), 63-177

\bibitem{Lamb3} J.S.W. Lamb and G.R.W. Quispel, ``Reversing k-symmetries in dynamical systems'', {\it Physica D} {\bf 73} (1994), 277-304

\bibitem{Lamb4} J.S.W. Lamb and G.R.W. Quispel, ``Cyclic reversing k-symmetry groups'', {\it Nonlinearity} {\bf  8} (1995), 1005-1026

\bibitem{QuiRob} G.R.W. Quispel and J.A.G. Roberts, ``Conservative and dissipative behaviour in reversible dynamical systems'', {\it Phys. Lett. A} {\bf 135} (1989), 337-342

\bibitem{LambNic} J.S.W. Lamb and M. Nicol, ``On symmetric attractors in reversible dynamical systems'', {\it  Physica D} {\bf 112} (1998), 281-297

\bibitem{BLH} H. Brands, J.S.W. Lamb and I. Hoveijn, ``Periodic orbits in k-symmetric dynamical systems'', {\it  Physica D} {\bf 84} (1995), 460-475

\bibitem{BPV} Biferale, L., D. Pierotti and A. Vulpiani, ``Time-reversible dynamical systems for turbulence'', {\it J. Phys. A} {\bf 31} (1998), 21-32

\bibitem{Dias} F. Dias and G. Iooss, ``Water-waves as a spatial dynamical system'', in {\it Handbook of Mathematical Fluid Dynamics. Vol. 2.} (pp 443-499), North-Holland, 2003






\bibitem{RodVan} H.M. Rodrigues and A. Vanderbauwhede, ``Symmetric perturbations of nonlinear equations: symmetry of small solutions'', {\it Nonlin. Anal. TMA} {\bf 2} (1978), 27-46

\bibitem{Mont} J.A. Montaldi, R.M.  Roberts, and I.N. Stewart. ``Periodic solutions near equilibria of symmetric Hamiltonian systems'', {\it Phil. Trans. Royal Soc. London A} {\bf 325} (1988), 237-293

\bibitem{CraKno} J.D. Crawford and E. Knobloch, ``Symmetry and symmetry-breaking bifurcations in fluid dynamics." Annual Review of Fluid Mechanics 23.1 (1991): 341-387.

\bibitem{Swift} J.W. Swift, ``Hopf bifurcation with the symmetry of the square'', {\it Nonlinearity} {\bf 1}  (1988), 333-377

\bibitem{Robbins} J.M. Robbins, ``Discrete symmetries in periodic-orbit theory'', {\it Phys. Rev. A} {\bf 40} (1989),  2128-2136

\bibitem{Zhil} B.I. Zhilinskiı, ``Symmetry, invariants, and topology in molecular models'', {\it Phys. Rep.} {\bf 341} (2001), 85-171

\bibitem{Ruck} A.M. Rucklidge and M. Silber, ``Bifurcations of periodic orbits with spatio-temporal symmetries'',  {\it Nonlinearity} {\bf 11} (1998), 1435-1455

\bibitem{VerGal} F. Verhulst, ``Discrete symmetric dynamical systems at the main resonances with application to axi-symmetric galaxies'', {\it Phil. Trans. Royal Soc. London A} {\bf 290} (1979), 435-465



\bibitem{Wal19} S. Walcher, ``Symmetries of ordinary differential equations: A short introduction'', {\tt arXiv:1911.01053} (2019)

\bibitem{KinWal} M.K. Kinyon, and S. Walcher, ``On ordinary differential equations admitting a finite linear group of symmetries'', {\it J. Math. Anal. Appl.} {\bf 216} (1997), 180-196

\bibitem{SchroWal} R. Schroeders, and S. Walcher, ``Orbit space reduction and localizations'', {\it Indag. Math.} {\bf 27} (2016), 1265-1278

\bibitem{GSW} F.D. Grosshans, J. Scheurle, and S. Walcher, ``Invariant sets forced by symmetry'', {\it J. Geom. Mech.} {\bf 4} (2012), 271-296

\bibitem{SchWal} J. Scheurle and S. Walcher, ``Hamiltonian Symmetry Reduction via Localizations: Theory and Application to a Barbell System'', {\it Acta Appl. Math.} {\bf 162} (2019), 121-143

\bibitem{HadWal} K.P. Hadeler, and Sebastian Walcher, ``Reducible ordinary differential equations'', {\it J. Nonlin. Sci.} {\bf 16} (2006), 583-613

\bibitem{CLPW} C. Christopher, J. Llibre, C. Pantazi, and S. Walcher, ``Inverse problems in Darboux theory of integrability'',{\it Acta Appl. Math} {\bf 120} (2012), 101-126

\bibitem{Goeke} A. Goeke, S. Walcher and E. Zerz, ``Classical quasi-steady state reduction. A mathematical characterization'', {\it Physica D} {\bf 345} (2017), 11-26

\bibitem{Broad1}  P. Broadbridge, B.H. Bradshaw-Hajek and D. Triadis, ``Exact non-classical symmetry solutions of Arrhenius reaction-diffusion'', {\it Proc. Royal Soc. A} {\bf 471} (2015), 20150580.

\bibitem{Broad2} P. Broadbridge, B.H. Bradshaw-Hajek and A.J. Hutchinson, ``Conditionally integrable PDEs, non-classical symmetries and applications'', {\it Proc. Royal Soc. A} {\bf 479} (2023), 20230209

\bibitem{Broad3} P. Broadbridge, R.M. Cherniha and J.M. Goard, ``Exact nonclassical symmetry solutions of Lotka-Volterra-type population systems'', {\it Eur. J. Appl. Math.} {\bf 34} (2023), 998-1016




\end{thebibliography}
\end{document}